\documentclass{article}
\usepackage{amsfonts,amsmath,amssymb,cite,graphicx,mathrsfs}
\numberwithin{equation}{section}
\begin{document}
\title{On the derivative of the associated Legendre function of the
first kind of integer order with respect to its degree (with
applications to the construction of the associated Legendre function
of the second kind of integer degree and order)}
\author{Rados{\l}aw Szmytkowski \\*[3ex]
Atomic Physics Division, 
Department of Atomic Physics and Luminescence, \\
Faculty of Applied Physics and Mathematics,
Gda{\'n}sk University of Technology, \\
Narutowicza 11/12, PL 80--233 Gda{\'n}sk, Poland \\
email: radek@mif.pg.gda.pl}
\date{\today}
\maketitle
\begin{abstract}
In our recent works [R.\ Szmytkowski, J.\ Phys.\ A 39 (2006) 15147;
corrigendum: 40 (2007) 7819; addendum: 40 (2007) 14887], we have
investigated the derivative of the Legendre function of the first
kind, $P_{\nu}(z)$, with respect to its degree $\nu$. In the present
work, we extend these studies and construct several representations
of the derivative of the associated Legendre function of the first
kind, $P_{\nu}^{\pm m}(z)$, with respect to the degree $\nu$, for
$m\in\mathbb{N}$. At first, we establish several contour-integral
representations of $\partial P_{\nu}^{\pm m}(z)/\partial\nu$. They
are then used to derive Rodrigues-type formulas for $[\partial
P_{\nu}^{\pm m}(z)/\partial\nu]_{\nu=n}$ with $n\in\mathbb{N}$. Next,
some closed-form expressions for $[\partial P_{\nu}^{\pm
m}(z)/\partial\nu]_{\nu=n}$ are obtained. These results are applied
to find several representations, both explicit and of the Rodrigues
type, for the associated Legendre function of the second kind of
integer degree and order, $Q_{n}^{\pm m}(z)$; the explicit
representations are suitable for use for numerical purposes in
various regions of the complex $z$-plane. Finally, the derivatives
$[\partial^{2}P_{\nu}^{m}(z)/\partial\nu^{2}]_{\nu=n}$,
$[\partial Q_{\nu}^{m}(z)/\partial\nu]_{\nu=n}$ and $[\partial
Q_{\nu}^{m}(z)/\partial\nu]_{\nu=-n-1}$, all with $m>n$, are
evaluated in terms of $[\partial P_{\nu}^{-m}(\pm 
z)/\partial\nu]_{\nu=n}$. \\*[1ex]
\end{abstract}
%
%
\section{Introduction}
\label{I}
\begin{figure}
\begin{center}
\includegraphics[bb=177 606 364 792]{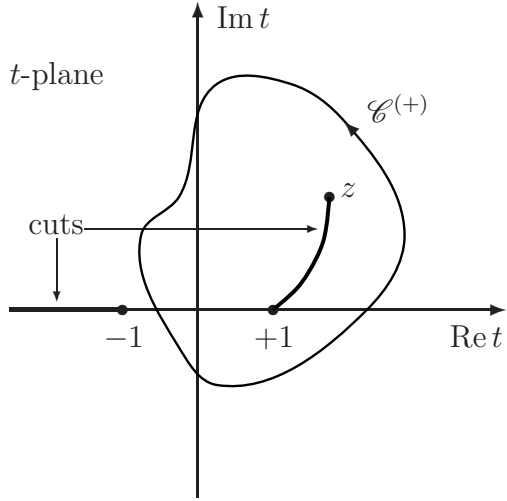}
\end{center}
\caption{The complex $t$-plane and the integration contour
$\mathscr{C}^{(+)}$ for Eq.\ (\ref{1.1}) and for the definition
(\ref{3.1}) of the associated Legendre function $P_{\nu}^{\pm m}(z)$
of the first kind of integer order. The cut joining the points $t=+1$
and $t=z$ (being two out of four branch points of the integrands in
Eqs.\ (\ref{1.1}) and (\ref{3.1})) is the circular arc (\ref{3.3}).}
\label{fig1}
\end{figure}

\setcounter{equation}{0}
In our recent paper \cite{Szmy06a}, we have proved that the
derivative of the Legendre function of the first kind, $P_{\nu}(z)$,
with respect to its degree $\nu$ may be given in the form
\begin{equation}
\frac{\partial P_{\nu}(z)}{\partial\nu}
=-P_{\nu}(z)\ln\frac{z+1}{2}+\frac{1}{2^{\nu}\pi\mathrm{i}}
\oint_{\mathscr{C}^{(+)}}\mathrm{d}t\:
\frac{(t^{2}-1)^{\nu}}{(t-z)^{\nu+1}}\ln\frac{t+1}{2},
\label{1.1}
\end{equation}
where the integration contour $\mathscr{C}^{(+)}$ is shown in Fig.\
\ref{fig1}.
Using the representation (\ref{1.1}), we have re-derived the
Rodrigues-type formula
\begin{equation}
\frac{\partial P_{\nu}(z)}{\partial\nu}\bigg|_{\nu=n}
=-P_{n}(z)\ln\frac{z+1}{2}
+\frac{1}{2^{n-1}n!}\frac{\mathrm{d}^{n}}{\mathrm{d}z^{n}}
\left[(z^{2}-1)^{n}\ln\frac{z+1}{2}\right]
\qquad (n\in\mathbb{N}),
\label{1.2}
\end{equation}
which was first obtained long ago by Jolliffe \cite{Joll19} with no
resort to the complex integration technique. Then we have shown that
$[\partial P_{\nu}(z)/\partial\nu]_{\nu=n}$ may be written as
\begin{equation}
\frac{\partial P_{\nu}(z)}{\partial\nu}\bigg|_{\nu=n}
=P_{n}(z)\ln\frac{z+1}{2}+R_{n}(z),
\label{1.3}
\end{equation}
where $P_{n}(z)$ is a Legendre polynomial and $R_{n}(z)$ is a
polynomial in $z$ of degree $n$. Using various techniques, we have
found several explicit representations of $R_{n}(z)$. In addition to
the representation
\begin{equation}
R_{n}(z)=2[\psi(2n+1)-\psi(n+1)]P_{n}(z)
+2\sum_{k=0}^{n-1}(-)^{k+n}\frac{2k+1}{(n-k)(k+n+1)}P_{k}(z),
\label{1.4}
\end{equation}
which is a slight modification of the one found by Bromwich
\cite{Brom13}, and the representation
\begin{equation}
R_{n}(z)=-2\psi(n+1)P_{n}(z)
+2\sum_{k=0}^{n}\frac{(k+n)!\psi(k+n+1)}{(k!)^{2}(n-k)!}
\left(\frac{z-1}{2}\right)^{k}
\label{1.5}
\end{equation}
due to Schelkunoff \cite{Sche41}, we have obtained the formula
\begin{equation}
R_{n}(z)=2\sum_{k=0}^{n}(-)^{k+n}\frac{(k+n)!}{(k!)^{2}(n-k)!}
[\psi(k+n+1)-\psi(k+1)]\left(\frac{z+1}{2}\right)^{k}.
\label{1.6}
\end{equation}
Here and henceforth, 
\begin{equation}
\psi(\zeta)=\frac{1}{\Gamma(\zeta)}
\frac{\mathrm{d}\Gamma(\zeta)}{\mathrm{d}\zeta}
\label{1.7}
\end{equation}
denotes the digamma function \cite{Davi65,Magn66,Grad94}. 

Still more recently, in the addendum \cite{Szmy07}, we have exploited
the Jolliffe's formula (\ref{1.2}) to present a derivation of Eq.\
(\ref{1.6}) being much simpler than the original one in
\cite{Szmy06a}. In addition, Eq.\ (\ref{1.2}) has been used to obtain
two further new representations of $R_{n}(z)$,
namely\footnote[1]{\label{FOOT1}~In \cite{Szmy07} the representations
(\ref{1.8}) and (\ref{1.9}) have been given in slightly different
forms.}
\begin{equation}
R_{n}(z)=2\psi(n+1)P_{n}(z)-2\left(\frac{z+1}{2}\right)^{n}
\sum_{k=0}^{n}{n \choose k}^{2}\psi(n-k+1)
\left(\frac{z-1}{z+1}\right)^{k}
\label{1.8}
\end{equation}
and
\begin{equation}
R_{n}(z)=2\psi(n+1)P_{n}(z)-2\left(\frac{z-1}{2}\right)^{n}
\sum_{k=0}^{n}{n \choose k}^{2}\psi(k+1)
\left(\frac{z+1}{z-1}\right)^{k}.
\label{1.9}
\end{equation}

In the course of solving boundary value problems of theoretical
acoustics, electromagnetism, heat conduction and some other branches
of theoretical physics and applied mathematics (cf, e.g.,
\cite{Cars13,Cars14a,Cars14b,Macd15,Cars21,Cars47,Smyt50,Fels55,
Bail56,Fels57,Fels59,Jone64,Bowm69,Fels73,Ariy85,Jone86,Baue87,Baue92,
Broa99,VanB07}), one occasionally encounters the derivative of the
associated Legendre function of the first kind
\cite{Todh75,Ferr77,Neum78,Hein78,Hein81,Olbr87,Byer93,Hobs96,Wang04,
Barn07,Wang21,Hobs31,Pras30,Pras32,Snow52,Lens54,Robi57,Robi58,Robi59,
Krat60,MacR67,Wang89,Temm96,Magn66,Grad94,Magn48,Erde53,Jahn60,Steg65,
Prud83,Prud03,Bryc08} of integer order with respect to its degree,
$\partial P_{\nu}^{\pm m}(z)/\partial\nu$ with $m\in\mathbb{N}$. We
have found that except for the case of $m=0$ studied in
\cite{Szmy06a,Szmy07}, thus far this derivative has not been a
subject of a systematic, exhaustive investigation and that the
relevant knowledge is very incomplete (cf section \ref{II}). It is
the purpose of the present paper to fill in this gap and to give an
extension of our previously-obtained results for $\partial
P_{\nu}(z)/\partial\nu$ to the case of $\partial P_{\nu}^{\pm
m}(z)/\partial\nu$. In particular, we shall extensively investigate
the functions $[\partial P_{\nu}^{\pm m}(z)/\partial\nu]_{\nu=n}$
with $n\in\mathbb{N}$. Besides being interesting for their own sake
and potentially useful for applications to the problems mentioned
above, the results will also allow us to contribute to the theory of
the associated Legendre function of the second kind of integer degree
and order, $Q_{n}^{\pm m}(z)$. In a simple manner, we shall obtain
several representations, both explicit and of the Rodrigues type, of
the latter function for various ranges of $n$ and $m$; we believe
that some of these representations are new. The explicit expressions
we provide for both $[\partial P_{\nu}^{\pm
m}(z)/\partial\nu]_{\nu=n}$ and $Q_{n}^{\pm m}(z)$ are suitable for
use for numerical purposes in various regions of the complex
$z$-plane.

The structure of the paper is as follows. Section \ref{II} provides a
summary of fragmentary research on $\partial
P_{\nu}^{\mu}(z)/\partial\nu$ done thus far by other authors. In
section \ref{III}, we make an overview of these properties of the
associated Legendre function of the first kind of integer order,
$P_{\nu}^{\pm m}(z)$, which will find applications in later parts of
the work. In section \ref{IV}, we find several contour-integral
representations of $\partial P_{\nu}^{\pm m}(z)/\partial\nu$. In
section \ref{V}, we investigate $[\partial P_{\nu}^{\pm
m}(z)/\partial\nu]_{\nu=n}$ with $n\in\mathbb{N}$, the cases
$0\leqslant m\leqslant n$ and $m>n$ being considered separately. In
each of these two cases, at first we use contour integrals for
$[\partial P_{\nu}^{\pm m}(z)/\partial\nu]_{\nu=n}$ to obtain
Rodrigues-type formulas for this function. Then, these formulas are
used to construct several closed-form representations of $[\partial
P_{\nu}^{\pm m}(z)/\partial\nu]_{\nu=n}$. Applications of the results
of section \ref{V} to the construction of some representations of the
associated Legendre function of the second kind of integer degree and
order constitute section \ref{VI.1}, while in sections \ref{VI.2} and
\ref{VI.3} the derivatives
$[\partial^{2}P_{\nu}^{m}(z)/\partial\nu^{2}]_{\nu=n}$, $[\partial
Q_{\nu}^{m}(z)/\partial\nu]_{\nu=n}$ and $[\partial
Q_{\nu}^{m}(z)/\partial\nu]_{\nu=-n-1}$, all with $m>n$, are
expressed in terms of $[\partial P_{\nu}^{-m}(\pm
z)/\partial\nu]_{\nu=n}$. The paper ends with an appendix, in which
some formulas for the Jacobi polynomials, exploited in sections
\ref{III} and \ref{V}, are listed.

Throughout the paper, we shall be adopting the standard convention,
according to which $z\in\mathbb{C}\setminus[-1,1]$, with the phases
restricted by
\begin{equation}
-\pi<\arg(z)<\pi, \qquad -\pi<\arg(z\pm1)<\pi
\label{1.10}
\end{equation}
(this corresponds to drawing a cut in the $z$-plane along the real
axis from $z=+1$ to $z=-\infty$), hence,
\begin{equation}
-z=\mathrm{e}^{\mp\mathrm{i}\pi}z
\qquad
-z+1=\mathrm{e}^{\mp\mathrm{i}\pi}(z-1)
\qquad
-z-1=\mathrm{e}^{\mp\mathrm{i}\pi}(z+1)
\qquad (\arg(z)\gtrless0).
\label{1.11}
\end{equation}
Also, it will be implicit that $x\in[-1,1]$, $\nu\in\mathbb{C}$,
$\mu\in\mathbb{C}$, $n\in\mathbb{N}$, $m\in\mathbb{N}$. Finally, it
will be understood that if the upper limit of a sum is less by unity
than the lower one, then the sum vanishes identically.

The definitions of the associated Legendre functions of the first and
second kinds used in the paper are those of Hobson \cite{Hobs31}.

The present paper may be considered as a complement to
\cite{Szmy09a}, where the derivative $\partial
P_{n}^{\mu}(z)/\partial\mu$ has been investigated exhaustively.
%
%
\section{Overview of research done on $\partial
P_{\nu}^{\mu}(z)/\partial\nu$ and $\partial
P_{\nu}^{\mu}(x)/\partial\nu$}
\label{II} 
\setcounter{equation}{0}
An overview of the research done on $\partial
P_{\nu}^{0}(z)/\partial\nu$ ($\equiv\partial P_{\nu}(z)/\partial\nu$)
was presented in \cite{Szmy06a}; to the references cited therein, one
should add the works
\cite{Fels55,Macd00,Macd02,Lebe39,Grin48,Jean60,Bowm69,Sakm70,
Nort71,vdPo87,Baue92,Soki97,From07,Coff08}. As regards the derivative
$\partial P_{\nu}^{\mu}(z)/\partial\nu$ with $\mu\neq0$, our
literature search showed very few results. Robin \cite{Robi56} (cf
also \cite[pp.\ 170--4]{Robi58}), differentiating term by term the
following series representation of $P_{\nu}^{\mu}(z)$
\cite{Magn66,Grad94}:
\begin{equation}
P_{\nu}^{\mu}(z)=\left(\frac{z+1}{z-1}\right)^{\mu/2}
\sum_{k=0}^{\infty}
\frac{\Gamma(k+\nu+1)}{k!\Gamma(\nu-k+1)\Gamma(k-\mu+1)}
\left(\frac{z-1}{2}\right)^{k}
\qquad (|z-1|<2)
\label{2.1}
\end{equation}
arrived at the formula\footnote[2]{\label{FOOT2}~Actually, Robin
\cite{Robi58,Robi56} used a definition of the digamma function
different from that in our Eq.\ (\ref{1.7}); his definition was
\begin{displaymath} 
\psi(\zeta)=\frac{1}{\Gamma(\zeta+1)}
\frac{\mathrm{d}\Gamma(\zeta+1)}{\mathrm{d}\zeta}. 
\end{displaymath}
Hence, our Eqs.\ (\ref{2.2}) and (\ref{2.6}) seemingly differ from
their counterparts in \cite{Robi58,Robi56}.}
\begin{eqnarray}
\frac{\partial P_{\nu}^{\mu}(z)}{\partial\nu}
&=& \left(\frac{z+1}{z-1}\right)^{\mu/2}\sum_{k=1}^{\infty}
\frac{\Gamma(k+\nu+1)}{k!\Gamma(\nu-k+1)\Gamma(k-\mu+1)}
\nonumber \\
&& \times[\psi(k+\nu+1)-\psi(\nu-k+1)]
\left(\frac{z-1}{2}\right)^{k}
\qquad (|z-1|<2).
\label{2.2}
\end{eqnarray}
If in Eq.\ (\ref{2.2}) one makes use of the well-known
\cite{Davi65,Magn66,Grad94} identities
\begin{equation}
\Gamma(\zeta)\Gamma(1-\zeta)=\frac{\pi}{\sin(\pi\zeta)}
\label{2.3}
\end{equation}
and
\begin{equation}
\psi(\zeta)=\psi(1-\zeta)-\pi\cot(\pi\zeta)
\label{2.4}
\end{equation}
and exploits Eq.\ (\ref{2.1}), this yields \cite[page 178]{Magn66}
\begin{eqnarray}
\frac{\partial P_{\nu}^{\mu}(z)}{\partial\nu}
&=& \pi\cot(\pi\nu)P_{\nu}^{\mu}(z)
-\frac{\sin(\pi\nu)}{\pi}\left(\frac{z+1}{z-1}\right)^{\mu/2}
\sum_{k=0}^{\infty}(-)^{k}
\frac{\Gamma(k+\nu+1)\Gamma(k-\nu)}{k!\Gamma(k-\mu+1)}
\nonumber \\
&& \times[\psi(k+\nu+1)-\psi(k-\nu)]\left(\frac{z-1}{2}\right)^{k}
\qquad (|z-1|<2).
\label{2.5}
\end{eqnarray}
Manipulating with the series on the right-hand side of the above
formula, Robin \cite{Robi58,Robi56} showed that for $\nu=n$ the
formula goes over into\footnote[3]{\label{FOOT3}~Equation (2) in
\cite{Robi56}, which is the counterpart of our Eq.\ (\ref{2.6}), was
misprinted: in front of the term containing the ${}_{3}F_{2}$
function, the factor $[(\mu+1)/(\mu-1)]^{m/2}$ (the original notation
of Robin is used here) is missing. In \cite[Eq.\ (329) on pp.\
171--2]{Robi58} the same formula was already printed correctly.}
\begin{eqnarray}
\frac{\partial P_{\nu}^{\mu}(z)}{\partial\nu}\bigg|_{\nu=n}
&=& \left(\frac{z+1}{z-1}\right)^{\mu/2}
\sum_{k=1}^{n}\frac{(k+n)!}{k!(n-k)!\Gamma(k-\mu+1)}
\nonumber \\
&& \quad \times[\psi(k+n+1)-\psi(n-k+1)]
\left(\frac{z-1}{2}\right)^{k}
\nonumber \\
&& +\,\frac{(2n+1)!}{(n+1)!\Gamma(n-\mu+2)}
\left(\frac{z+1}{z-1}\right)^{\mu/2}\left(\frac{z-1}{2}\right)^{n+1}
\nonumber \\
&& \quad \times{}_{3}F_{2}
\left(
\begin{array}{c}
1, 1, 2n+2 \\
n+2, n+2-\mu
\end{array}
;\frac{1-z}{2}
\right),
\label{2.6}
\end{eqnarray}
with ${}_{3}F_{2}$ being a generalized hypergeometric
function\footnote[4]{\label{FOOT4}~If in Eq.\ (\ref{2.6}) one sets
$\mu=0$, combines the result with the Schelkunoff's formula for
$[\partial P_{\nu}(z)/\partial\nu]_{\nu=n}$ following from Eqs.\
(\ref{1.3}) and (\ref{1.5}), then solves the emerging equation for
$_{3}F_{2}(1,1,2n+2;n+2,n+2;(1-z)/2)$ and replaces therein $z$ by
$1-2z$, one obtains
\begin{eqnarray*}
_{3}F_{2}
\left(
\begin{array}{c}
1, 1, 2n+2 \\
n+2, n+2
\end{array}
;z
\right)
&=& (-)^{n+1}\frac{[(n+1)!]^{2}}{(2n+1)!z^{n+1}}
\bigg\{P_{n}(1-2z)\ln(1-z)
\nonumber \\
&& +\,\sum_{k=1}^{n}(-)^{k}\frac{(k+n)!}{(k!)^{2}(n-k)!}
[\psi(k+n+1)+\psi(n-k+1)-2\psi(n+1)]z^{k}\bigg\}.
\end{eqnarray*}
This relationship should replace Eq.\ (7.4.1.35) in \cite[pp.\
421--2]{Prud03}, which is incorrect in view of the fact that in
\cite[p.\ 685]{Prud03} the digamma function has been defined as in
our Eq.\ (\ref{1.7}) and \emph{not\/} as in
\cite{Sche41,Robi58,Robi56}.}. In the particular case of $n=0$, the
finite sum on the right-hand side of Eq.\ (\ref{2.6}) vanishes,
while the $_{3}F_{2}$ function appearing therein reduces to
$_{2}F_{1}(1,1;2-\mu;(1-z)/2)$, so that the equation goes over
into\footnote[5]{\label{FOOT5}~It seems worthwhile to add at this
place that Eq.\ (\ref{2.7}) may be used to express the closed-form
momentum-space representation of the nonrelativistic Coulomb Green
function, found by Hostler in \cite{Host64}, in terms of the
derivative $[\partial P_{\nu}^{\mu}(z)/\partial\nu]_{\nu=0}$ with
suitably chosen $\mu$ and $z$.} \cite[p.\ 177]{Magn66}
\begin{equation}
\frac{\partial P_{\nu}^{\mu}(z)}{\partial\nu}\bigg|_{\nu=0}
=\frac{z-1}{2\Gamma(2-\mu)}\left(\frac{z+1}{z-1}\right)^{\mu/2}
{}_{2}F_{1}\left(1,1;2-\mu;\frac{1-z}{2}\right).
\label{2.7}
\end{equation}
For the sake of completeness, we mention here that a formula for
$[\partial P_{\nu}^{\mu}(z)/\partial\nu]_{\nu=n-1/2}$, which we do
not display here, has been provided by Brychkov in \cite{Bryc08}.

Results on $\partial P_{\nu}^{\mu}(x)/\partial\nu$ are also scarce
and fragmentary. A formula analogous to Eq.\ (\ref{2.2}) may be found
in \cite[p.\ 1026]{Grad94} and \cite[p.\ 94]{Magn48}. Counterparts of
Eqs.\ (\ref{2.5}) and (\ref{2.7}) are given in \cite{Magn66} on pp.\
178 and 177, respectively. Closed-form expressions for $[\partial
P_{\nu}^{-1}(x)/\partial\nu]_{\nu=0}$ and
$[\partial P_{\nu}^{-1}(x)/\partial\nu]_{\nu=1}$ are presented in
\cite[pp.\ 1026--7]{Grad94}, \cite[p.\ 94]{Magn48} and \cite[p.\
335]{Steg65}. Tsu \cite{Tsu61} found explicit representations
of $[\partial P_{\nu}^{0}(x)/\partial\nu]_{\nu=0}$ ($\equiv[\partial
P_{\nu}(x)/\partial\nu]_{\nu=0}$) and $[\partial
P_{\nu}^{1}(x)/\partial\nu]_{\nu=0}$ and provided the recursive
relation
\begin{eqnarray}
(1-x^{2})^{1/2}\frac{\partial P_{\nu}^{m+1}(x)}
{\partial\nu}\bigg|_{\nu=n}
-(n-m)x\frac{\partial P_{\nu}^{m}(x)}{\partial\nu}\bigg|_{\nu=n}
+(n+m)\frac{\partial P_{\nu}^{m}(x)}{\partial\nu}\bigg|_{\nu=n-1}
\nonumber \\
=xP_{n}^{m}(x)-P_{n-1}^{m}(x),
\label{2.8}
\end{eqnarray}
enabling one to generate $[\partial
P_{\nu}^{m}(x)/\partial\nu]_{\nu=n}$ for other values of $n$ and
$m$; however, no general formula for $[\partial
P_{\nu}^{m}(x)/\partial\nu]_{\nu=n}$ was given in that work. Finally,
in a study on the Dirichlet averages of $x^{t}\ln x$, Carlson
\cite{Carl87} arrived at the following closed-form
representation\footnote[6]{\label{FOOT6}~We have transformed
Carlson's original formulas so that Eqs.\ (\ref{2.9}) and
(\ref{2.11}) are concurrent with the notation used in the rest of the
present paper.} of $[\partial P_{\nu}^{-m}(x)/\partial\nu]_{\nu=n}$
with $0\leqslant m\leqslant n$:
\begin{eqnarray}
\frac{\partial P_{\nu}^{-m}(x)}{\partial\nu}\bigg|_{\nu=n}
&=& P_{n}^{-m}(x)\ln\frac{1+x}{2}
-[\psi(n+m+1)+\psi(n+1)]P_{n}^{-m}(x)
\nonumber \\
&& +\,\frac{(n-m)!}{(n+m)!}\left(\frac{1-x^{2}}{4}\right)^{m/2}
\nonumber \\
&& \quad\times\sum_{k=0}^{n-m}(-)^{k}\frac{(k+n+m)!\psi(k+n+m+1)}
{k!(k+m)!(n-m-k)!}\left(\frac{1-x}{2}\right)^{k}
\nonumber \\
&& +\,\left(\frac{1-x}{1+x}\right)^{m/2}
\sum_{k=0}^{n}(-)^{k}\frac{(k+n)!\psi(k+n+1)}{k!(k+m)!(n-k)!}
\left(\frac{1-x}{2}\right)^{k}
\qquad (0\leqslant m\leqslant n)
\nonumber \\
&&
\label{2.9}
\end{eqnarray}
and proved the identity
\begin{eqnarray}
\frac{\partial P_{\nu}^{m}(x)}{\partial\nu}\bigg|_{\nu=n}
&=& (-)^{m}\frac{(n+m)!}{(n-m)!}
\frac{\partial P_{\nu}^{-m}(x)}{\partial\nu}\bigg|_{\nu=n}
\nonumber \\
&& +\,[\psi(n+m+1)-\psi(n-m+1)]P_{n}^{m}(x)
\qquad (0\leqslant m\leqslant n);
\label{2.10}
\end{eqnarray}
in addition, he showed that
\begin{eqnarray}
\frac{\partial P_{\nu}^{-m}(x)}{\partial\nu}\bigg|_{\nu=0}
&=& (-)^{m}P_{0}^{-m}(-x)\ln\frac{1+x}{2}
\nonumber \\
&& +\,\psi(1)P_{0}^{-m}(x)-(-)^{m}\psi(m+1)P_{0}^{-m}(-x)
\nonumber  \\
&& +\,(-)^{m}\left(\frac{1-x^{2}}{4}\right)^{-m/2}
\sum_{k=0}^{m-1}(-)^{k}\frac{\psi(m-k+1)}{k!(m-k)!}
\left(\frac{1-x}{2}\right)^{k}
\qquad (m>0).
\nonumber \\
&&
\label{2.11}
\end{eqnarray}
%
%
\section{Definition and some relevant properties of the associated
Legendre function of the first kind of integer order}
\label{III}
\setcounter{equation}{0}
In this section, we shall present these properties of the associated
Legendre function of the first kind of integer order which will find
applications in later parts of this paper.
\subsection{Function of arbitrary degree}
\label{III.1}
The associated Legendre function of the first kind of complex degree
$\nu$ and integer order $\pm m$ may be defined as the following
generalization of the Schl{\"a}fli contour integral 
\cite[p.\ 191]{Hobs31}:
\begin{equation}
P_{\nu}^{\pm m}(z)=\frac{\Gamma(\nu\pm m+1)}{\Gamma(\nu+1)}
\frac{1}{2^{\nu+1}\pi\mathrm{i}}(z^{2}-1)^{\pm m/2}
\oint_{\mathscr{C}^{(+)}}\mathrm{d}t\:
\frac{(t^{2}-1)^{\nu}}{(t-z)^{\nu\pm m+1}}.
\label{3.1}
\end{equation}
The integration path $\mathscr{C}^{(+)}$, shown in Fig.\ \ref{fig1},
is a closed circuit enclosing the points $t=+1$ and $t=z$ in the
counter-clockwise sense. If $\nu$ is not an integer, the integrand
has four branch points located at $t=\pm1$, $t=z$ and
$|t|=\infty$. To make the integrand single-valued, we make two cuts
in the $t$-plane. The first one is
\begin{equation}
t_{\mathrm{cut}}=-\eta
\qquad (1\leqslant\eta<\infty),
\label{3.2}
\end{equation}
which is the semi-line drawn along the negative real semi-axis from
$t=-1$ to $t=-\infty$. The second one is the curve
\begin{equation}
t_{\mathrm{cut}}=\frac{1+\eta z}{z+\eta}
\qquad (1\leqslant\eta<\infty)
\label{3.3}
\end{equation}
joining the points $t=+1$ and $t=z$. It is seen to be that out of two
circular arcs of radius
\begin{equation}
\rho=\sqrt{1+\left(\frac{|z|^{2}-1}{2\textrm{Im}(z)}\right)^{2}},
\label{3.4}
\end{equation}
centered at 
\begin{equation}
t_{0}=\mathrm{i}\,\frac{|z|^{2}-1}{2\textrm{Im}(z)}
\label{3.5}
\end{equation}
and connecting the points $t=+1$ and $t=z$, which does \emph{not\/}
go through the point $t=-1$. The contour $\mathscr{C}^{(+)}$ is not
to cross any of the two cuts (\ref{3.2}) and (\ref{3.3}). The phases
in the integrand in Eq.\ (\ref{3.1}) are stipulated as follows: at
the point on the right to $t=+1$ (and on the right to $z$ if the
latter be real), where the path $\mathscr{C}^{(+)}$ crosses the real
axis, we set $\arg(t\pm1)=0$ and $|\arg(t-z)|<\pi$. In the plane with
the cross-cut along the real axis from $z=+1$ to $z=-\infty$ (cf the
remark below Eq.\ (\ref{1.10})), the function $P_{\nu}^{\pm m}(z)$ is
single-valued.

It is seen from Eq.\ (\ref{3.1}) that
\begin{equation}
P_{\nu}^{m}(z)=(z^{2}-1)^{m/2}
\frac{\mathrm{d}^{m}P_{\nu}(z)}{\mathrm{d}z^{m}},
\label{3.6}
\end{equation}
where $P_{\nu}(z)\equiv P_{\nu}^{0}(z)$ is the Legendre function of
the first kind.

It may be shown \cite{Hobs31} that $P_{\nu}^{\pm m}(z)$ possesses the
property
\begin{equation}
P_{-\nu-1}^{\pm m}(z)=P_{\nu}^{\pm m}(z).
\label{3.7}
\end{equation}
Replacing in Eq.\ (\ref{3.1}) $\nu$ with $-\nu-1$, exploiting the
fact that
\begin{equation}
\frac{\Gamma(-\nu\pm m)}{\Gamma(-\nu)}
=(-)^{m}\frac{\Gamma(\nu+1)}{\Gamma(\nu\mp m+1)}
\label{3.8}
\end{equation}
and making use of Eq.\ (\ref{3.7}), yields
\begin{equation}
P_{\nu}^{\pm m}(z)=(-)^{m}\frac{\Gamma(\nu+1)}{\Gamma(\nu\mp m+1)}
\frac{2^{\nu}}{\pi\mathrm{i}}(z^{2}-1)^{\pm m/2}
\oint_{\mathscr{C}^{(+)}}\mathrm{d}t\:
\frac{(t-z)^{\nu\mp m}}{(t^{2}-1)^{\nu+1}}.
\label{3.9}
\end{equation}
If in Eq.\ (\ref{3.9}) we change the integration variable from $t$ to
\begin{equation}
u=-1+2\frac{z+1}{t+1},
\label{3.10}
\end{equation}
this results in
\begin{equation}
P_{\nu}^{\pm m}(z)=\frac{\Gamma(\nu+1)}{\Gamma(\nu\mp m+1)}
\frac{1}{2^{\nu+1}\pi\mathrm{i}}
\left(\frac{z-1}{z+1}\right)^{\pm m/2}
\oint_{\mathscr{C}^{\prime\,(+)}}\mathrm{d}u\:
\frac{(u-1)^{\nu\mp m}(u+1)^{\nu\pm m}}{(u-z)^{\nu+1}}.
\label{3.11}
\end{equation}
The contour $\mathscr{C}^{\prime\,(+)}$ surrounds the points $u=+1$
and $u=z$ in the counter-clockwise sense and does not cross either of
the two straight-line cuts
\begin{equation}
u_{\mathrm{cut}}=-1+2\frac{z+1}{1-\eta}
\qquad (1\leqslant\eta<\infty)
\label{3.12}
\end{equation}
and
\begin{equation}
u_{\mathrm{cut}}=1+2\frac{z-1}{1+\eta}
\qquad (1\leqslant\eta<\infty).
\label{3.13}
\end{equation}
Likewise, if in Eq.\ (\ref{3.9}) the variable $t$ is replaced by
\begin{equation}
u=1-2\frac{z-1}{t-1},
\label{3.14}
\end{equation}
one finds
\begin{equation}
P_{\nu}^{\pm m}(z)=\frac{\Gamma(\nu+1)}{\Gamma(\nu\mp m+1)}
\frac{1}{2^{\nu+1}\pi\mathrm{i}}
\left(\frac{z+1}{z-1}\right)^{\pm m/2}
\oint_{\mathscr{C}^{\prime\prime\,(+)}}\mathrm{d}u\:
\frac{(u-1)^{\nu\pm m}(u+1)^{\nu\mp m}}{(u-z)^{\nu+1}},
\label{3.15}
\end{equation}
where the path $\mathscr{C}^{\prime\prime\,(+)}$ encloses the points
$u=+1$ and $u=z$ in the counter-clockwise sense and also does not
cross the cuts (\ref{3.12}) and (\ref{3.13}). It is evident that the
contour $\mathscr{C}^{\prime\prime\,(+)}$ may be deformed into the
contour $\mathscr{C}^{\prime\,(+)}$ without changing the value of the
integral in Eq.\ (\ref{3.15}), i.e., it holds that
\begin{equation}
P_{\nu}^{\pm m}(z)=\frac{\Gamma(\nu+1)}{\Gamma(\nu\mp m+1)}
\frac{1}{2^{\nu+1}\pi\mathrm{i}}
\left(\frac{z+1}{z-1}\right)^{\pm m/2}
\oint_{\mathscr{C}^{\prime\,(+)}}\mathrm{d}u\:
\frac{(u-1)^{\nu\pm m}(u+1)^{\nu\mp m}}{(u-z)^{\nu+1}}.
\label{3.16}
\end{equation}

As a corollary, from Eqs.\ (\ref{3.11}) and (\ref{3.16}) one obtains
the relation
\begin{equation}
P_{\nu}^{-m}(z)
=\frac{\Gamma(\nu-m+1)}{\Gamma(\nu+m+1)}P_{\nu}^{m}(z).
\label{3.17}
\end{equation}
If this is combined with Eq.\ (\ref{3.9}), this results in
\begin{equation}
P_{\nu}^{\pm m}(z)=(-)^{m}\frac{\Gamma(\nu+1)}{\Gamma(\nu\mp m+1)}
\frac{2^{\nu}}{\pi\mathrm{i}}(z^{2}-1)^{\mp m/2}
\oint_{\mathscr{C}^{(+)}}\mathrm{d}t\:
\frac{(t-z)^{\nu\pm m}}{(t^{2}-1)^{\nu+1}}.
\label{3.18}
\end{equation}
(It is worthwhile to add that Eq.\ (\ref{3.18}) may be also obtained
from equation (\ref{3.9}) by subjecting the integral in the latter to
the variable transformation (\ref{4.2}) and deforming suitably the
resulting integration contour. Then Eq.\ (\ref{3.17}) appears to be a
corollary from Eqs.\ (\ref{3.9}) and (\ref{3.18}).)

On the cut $-1\leqslant x\leqslant+1$, after Hobson \cite{Hobs31}, it
is customary to define
\begin{eqnarray}
P_{\nu}^{\pm m}(x) &=& \mathrm{e}^{\pm\mathrm{i}\pi m/2}
P_{\nu}^{\pm m}(x+\mathrm{i}0)
=\mathrm{e}^{\mp\mathrm{i}\pi m/2}P_{\nu}^{\pm m}(x-\mathrm{i}0)
\nonumber \\
&=& \frac{1}{2}\left[\mathrm{e}^{\pm\mathrm{i}\pi m/2}
P_{\nu}^{\pm m}(x+\mathrm{i}0)+\mathrm{e}^{\mp\mathrm{i}\pi m/2}
P_{\nu}^{\pm m}(x-\mathrm{i}0)\right].
\label{3.19}
\end{eqnarray}
%
%
\subsection{Function of integer degree}
\label{III.2}
If $\nu=n$, the cut (\ref{3.3}) in the definition of the contour
integral (\ref{3.1}) is unnecessary and may be removed. Then, by the
theory of residues, from Eq.\ (\ref{3.1}) one has the Rodrigues-type
formula
\begin{equation}
P_{n}^{\pm m}(z)=\frac{1}{2^{n}n!}(z^{2}-1)^{\pm m/2}
\frac{\mathrm{d}^{n\pm m}}{\mathrm{d}z^{n\pm m}}(z^{2}-1)^{n},
\label{3.20}
\end{equation}
subject to the constraint $0\leqslant m\leqslant n$ if the lower
signs are chosen. If Eq.\ (\ref{3.20}) is combined with
\begin{equation}
P_{n}^{-m}(z)=\frac{(n-m)!}{(n+m)!}P_{n}^{m}(z)
\qquad (0\leqslant m\leqslant n),
\label{3.21}
\end{equation}
which is the direct consequence of Eq.\ (\ref{3.17}), this yields
\begin{equation}
P_{n}^{\pm m}(z)=\frac{1}{2^{n}n!}
\frac{(n\pm m)!}{(n\mp m)!}(z^{2}-1)^{\mp m/2}
\frac{\mathrm{d}^{n\mp m}}{\mathrm{d}z^{n\mp m}}(z^{2}-1)^{n}
\qquad (0\leqslant m\leqslant n).
\label{3.22}
\end{equation}

For $m>n$, Eq.\ (\ref{3.20}) implies
\begin{equation}
P_{n}^{m}(z)=0
\qquad (m>n).
\label{3.23}
\end{equation}
Another property of $P_{n}^{\pm m}(z)$, which shall prove to be
useful in later considerations, is
\begin{equation}
P_{n}^{\pm m}(-z)=(-)^{n}P_{n}^{\pm m}(z)
\qquad (0\leqslant m\leqslant n).
\label{3.24}
\end{equation}
This may be obtained from Eq.\ (\ref{3.20}), with the aid of the
relations (\ref{1.11}).

Other Rodrigues-type representations of $P_{n}^{\pm m}(z)$ may be
obtained by applying the theory of residues to the contour integrals
(\ref{3.11}) and (\ref{3.16}), after setting therein $\nu=n$ and
removing, now redundant, the cut (\ref{3.13}). For $0\leqslant
m\leqslant n$, this renders two formulas
\begin{equation}
P_{n}^{\pm m}(z)=\frac{1}{2^{n}(n\mp m)!}
\left(\frac{z-1}{z+1}\right)^{\pm m/2}
\frac{\mathrm{d}^{n}}{\mathrm{d}z^{n}}
\left[(z-1)^{n\mp m}(z+1)^{n\pm m}\right]
\qquad (0\leqslant m\leqslant n)
\label{3.25}
\end{equation}
and
\begin{equation}
P_{n}^{\pm m}(z)=\frac{1}{2^{n}(n\mp m)!}
\left(\frac{z+1}{z-1}\right)^{\pm m/2}
\frac{\mathrm{d}^{n}}{\mathrm{d}z^{n}}
\left[(z-1)^{n\pm m}(z+1)^{n\mp m}\right]
\qquad (0\leqslant m\leqslant n),
\label{3.26}
\end{equation}
obtained originally by Schendel \cite{Sche77} in a different way. If
$m>n$, proceeding in the analogous way, from Eq.\ (\ref{3.11}) one
finds
\begin{equation}
P_{n}^{-m}(z)=\frac{1}{2^{n}(n+m)!}
\left(\frac{z+1}{z-1}\right)^{m/2}
\frac{\mathrm{d}^{n}}{\mathrm{d}z^{n}}
\left[(z-1)^{n+m}(z+1)^{n-m}\right]
\qquad (m>n).
\label{3.27}
\end{equation}

For the sake of later applications, we shall derive here still
another Rodrigues-type representation of $P_{n}^{-m}(z)$ valid for
$m>n$. To this end, at first we observe that if the lower signs are
chosen in Eq.\ (\ref{3.18}) and if one sets $\nu=n$, the integrand in
the resulting equation becomes single-valued in the domain enclosed
by the contour $\mathscr{C}^{(+)}$ and is seen to possess two poles
in this region: one of order $n+1$, located at $t=+1$, and the other,
of order $m-n$, located at $t=z$. Thus, we may remove the cut
(\ref{3.3}) and, by the Cauchy theorem, write
\begin{eqnarray}
P_{n}^{-m}(z) &=& 
(-)^{m}\frac{n!}{(n+m)!}\frac{2^{n}}{\pi\mathrm{i}}
(z^{2}-1)^{m/2}\oint_{\mathscr{C}_{+1}^{(+)}}\mathrm{d}t\:
\frac{(t-z)^{n-m}}{(t^{2}-1)^{n+1}}
\nonumber \\
&& +\,(-)^{m}\frac{n!}{(n+m)!}\frac{2^{n}}{\pi\mathrm{i}}
(z^{2}-1)^{m/2}\oint_{\mathscr{C}_{z}^{(+)}}\mathrm{d}t\:
\frac{(t-z)^{n-m}}{(t^{2}-1)^{n+1}}
\qquad (m>n).
\label{3.28}
\end{eqnarray}
Here, the path $\mathscr{C}_{+1}^{(+)}$ surrounds the point $t=+1$ in
the positive direction and leaves the points $t=-1$ and $t=z$
outside, while the contour $\mathscr{C}_{z}^{(+)}$ encloses the point
$t=z$, is run in the positive sense, with the points $t=\pm1$ located
exterior to it; none of the two paths crosses the cut (\ref{3.2}).
Instead of applying the theory of residues already at this stage, we
subject the first integral in Eq.\ (\ref{3.28}) to the variable
transformation (\ref{3.10}). This gives
\begin{eqnarray}
P_{n}^{-m}(z) &=& 
\frac{n!}{(n+m)!}\frac{1}{2^{n+1}\pi\mathrm{i}}
\left(\frac{z-1}{z+1}\right)^{m/2}
\oint_{\mathscr{C}_{z}^{\prime\,(+)}}\mathrm{d}u\:
\frac{(u-1)^{n-m}(u+1)^{n+m}}{(u-z)^{n+1}}
\nonumber \\
&& +\,(-)^{m}\frac{n!}{(n+m)!}\frac{2^{n}}{\pi\mathrm{i}}
(z^{2}-1)^{m/2}\oint_{\mathscr{C}_{z}^{(+)}}\mathrm{d}t\:
\frac{(t-z)^{n-m}}{(t^{2}-1)^{n+1}}
\qquad (m>n),
\label{3.29}
\end{eqnarray}
where the contour $\mathscr{C}_{z}^{\prime\,(+)}$ encircles the point
$u=z$ counter-clockwise, does not enclose the points $u=\pm1$ and
also does not cross the cut (\ref{3.12}). Applying now the theory of
residues to Eq.\ (\ref{3.29}), we obtain
\begin{eqnarray}
P_{n}^{-m}(z) &=& \frac{1}{2^{n}(n+m)!}
\left(\frac{z-1}{z+1}\right)^{m/2}
\frac{\mathrm{d}^{n}}{\mathrm{d}z^{n}}
\left[(z-1)^{n-m}(z+1)^{n+m}\right]
\nonumber \\
&& +\,(-)^{m}\frac{2^{n+1}n!}{(n+m)!(m-n-1)!}(z^{2}-1)^{m/2}
\frac{\mathrm{d}^{m-n-1}}{\mathrm{d}z^{m-n-1}}(z^{2}-1)^{-n-1}
\qquad (m>n),
\nonumber \\
&&
\label{3.30}
\end{eqnarray}
which is the desired result. Comparison with Eq.\ (\ref{3.27}) shows
that the first term on the right-hand side of Eq.\ (\ref{3.30})
equals $(-)^{n}P_{n}^{-m}(-z)$ (with $m>n$), and consequently
\begin{eqnarray}
P_{n}^{-m}(z) &=& (-)^{n}P_{n}^{-m}(-z)
\nonumber \\
&& +\,(-)^{m}\frac{2^{n+1}n!}{(n+m)!(m-n-1)!}(z^{2}-1)^{m/2}
\frac{\mathrm{d}^{m-n-1}}{\mathrm{d}z^{m-n-1}}(z^{2}-1)^{-n-1}
\qquad (m>n).
\nonumber \\
&&
\label{3.31}
\end{eqnarray}

Combining Eqs.\ (\ref{3.20}), (\ref{3.22}), (\ref{3.25}) to
(\ref{3.27}) and (\ref{3.31}) with Eq.\ (\ref{A.1}) and using,
whenever necessary, Eqs.\ (\ref{3.21}) and (\ref{3.24}), yields the
following formulas relating the Legendre function considered here to
particular Jacobi polynomials:
\begin{equation}
P_{n}^{m}(z)=\frac{n!}{(n-m)!}\left(\frac{z+1}{z-1}\right)^{m/2}
P_{n}^{(-m,m)}(z)
\qquad (0\leqslant m\leqslant n),
\label{3.32}
\end{equation}
\begin{equation}
P_{n}^{-m}(-z)=(-)^{n}\frac{n!}{(n+m)!}
\left(\frac{z+1}{z-1}\right)^{m/2}P_{n}^{(-m,m)}(z),
\label{3.33}
\end{equation}
\begin{equation}
P_{n}^{m}(z)=\frac{n!}{(n-m)!}\left(\frac{z-1}{z+1}\right)^{m/2}
P_{n}^{(m,-m)}(z)
\qquad (0\leqslant m\leqslant n),
\label{3.34}
\end{equation}
\begin{equation}
P_{n}^{-m}(z)=\frac{n!}{(n+m)!}\left(\frac{z-1}{z+1}\right)^{m/2}
P_{n}^{(m,-m)}(z),
\label{3.35}
\end{equation}
\begin{equation}
P_{n}^{m}(z)=\frac{(n+m)!}{n!}
\left(\frac{z^{2}-1}{4}\right)^{-m/2}P_{n+m}^{(-m,-m)}(z)
\qquad (0\leqslant m\leqslant n),
\label{3.36}
\end{equation}
\begin{equation}
P_{n}^{m}(z)=\frac{(n+m)!}{n!}
\left(\frac{z^{2}-1}{4}\right)^{m/2}P_{n-m}^{(m,m)}(z)
\qquad (0\leqslant m\leqslant n),
\label{3.37}
\end{equation}
\begin{equation}
P_{n}^{-m}(z)-(-)^{n}P_{n}^{-m}(-z)
=(-)^{m}\frac{n!}{(n+m)!}\left(\frac{z^{2}-1}{4}\right)^{-m/2}
P_{m-n-1}^{(-m,-m)}(z)
\qquad (m>n).
\label{3.38}
\end{equation}
We shall make extensive use of the above formulas in sections
\ref{V.1.2} and \ref{V.4.2}.
%
%
\section{Contour-integral representations of $\partial P_{\nu}^{\pm
m}(z)/\partial\nu$}
\label{IV}
\setcounter{equation}{0}
We begin our investigations on the derivative $\partial P_{\nu}^{\pm
m}(z)/\partial\nu$ with the derivation of its several
contour-integral representations.

Differentiation of Eq.\ (\ref{3.1}) with respect to $\nu$ gives the
first such representation:
\begin{eqnarray}
\frac{\partial P_{\nu}^{\pm m}(z)}{\partial\nu}
&=& [\psi(\nu\pm m+1)-\psi(\nu+1)]P_{\nu}^{\pm m}(z)
\nonumber \\
&& +\,\frac{\Gamma(\nu\pm m+1)}{\Gamma(\nu+1)}
\frac{1}{2^{\nu+1}\pi\mathrm{i}}(z^{2}-1)^{\pm m/2}
\oint_{\mathscr{C}^{(+)}}\mathrm{d}t\:
\frac{(t^{2}-1)^{\nu}}{(t-z)^{\nu\pm m+1}}
\ln\frac{t^{2}-1}{2(t-z)}.
\nonumber \\
&&
\label{4.1}
\end{eqnarray}

Consider now the following linear fractional transformation:
\begin{equation}
s=-1+(z+1)\frac{t-1}{t-z}.
\label{4.2}
\end{equation}
It maps the complex $t$-plane onto the complex $s$-plane. In
particular, the points $t=-1$, $+1$, $z$ and $\infty$ are mapped onto
the points $s=+1$, $-1$, $\infty$ and $z$, respectively, the cut
(\ref{3.2}) is mapped onto the cut
\begin{equation}
s_{\mathrm{cut}}=\frac{1+\eta z}{z+\eta}
\qquad (1\leqslant\eta<\infty)
\label{4.3}
\end{equation}
and the cut (\ref{3.3}) onto the cut
\begin{equation}
s_{\mathrm{cut}}=-\eta
\qquad (1\leqslant\eta<\infty).
\label{4.4}
\end{equation}
In addition, the path $\mathscr{C}^{(+)}$ is mapped onto the contour
$\mathscr{C}^{\prime\prime\prime\,(-)}$, which encloses the points
$s=1$ and $s=z$ in the \emph{clock}-wise sense and does not cross any
of the two cuts (\ref{4.3}) and (\ref{4.4}).

To obtain the second representation of $\partial P_{\nu}^{\pm
m}(z)/\partial\nu$, we rewrite Eq.\ (\ref{4.1}) in the form
\begin{eqnarray}
\frac{\partial P_{\nu}^{\pm m}(z)}{\partial\nu}
&=& [\psi(\nu\pm m+1)-\psi(\nu+1)]P_{\nu}^{\pm m}(z)
\nonumber \\
&& +\,\frac{\Gamma(\nu\pm m+1)}{\Gamma(\nu+1)}
\frac{1}{2^{\nu+1}\pi\mathrm{i}}(z^{2}-1)^{\pm m/2}
\oint_{\mathscr{C}^{(+)}}\mathrm{d}t\:
\frac{(t^{2}-1)^{\nu}}{(t-z)^{\nu\pm m+1}}\ln\frac{t+1}{2}
\nonumber \\
&& +\frac{\Gamma(\nu\pm m+1)}{\Gamma(\nu+1)}
\frac{1}{2^{\nu+1}\pi\mathrm{i}}(z^{2}-1)^{\pm m/2}
\oint_{\mathscr{C}^{(+)}}\mathrm{d}t\:
\frac{(t^{2}-1)^{\nu}}{(t-z)^{\nu\pm m+1}}\ln\frac{t-1}{t-z}.
\label{4.5}
\end{eqnarray}
Let us look closer at the second integral on the right-hand side of
Eq.\ (\ref{4.5}), which is
\begin{equation}
I_{\nu}^{\pm m}(z)=\frac{\Gamma(\nu\pm m+1)}{\Gamma(\nu+1)}
\frac{1}{2^{\nu+1}\pi\mathrm{i}}(z^{2}-1)^{\pm m/2}
\oint_{\mathscr{C}^{(+)}}\mathrm{d}t\:
\frac{(t^{2}-1)^{\nu}}{(t-z)^{\nu\pm m+1}}\ln\frac{t-1}{t-z}.
\label{4.6}
\end{equation}
Subjecting this integral to the transformation (\ref{4.2}) results in
\begin{equation}
I_{\nu}^{\pm m}(z)=-\frac{\Gamma(\nu\pm m+1)}{\Gamma(\nu+1)}
\frac{1}{2^{\nu+1}\pi\mathrm{i}}(z^{2}-1)^{\mp m/2}
\oint_{\mathscr{C}^{\prime\prime\prime\,(-)}}\mathrm{d}s\:
\frac{(s^{2}-1)^{\nu}}{(s-z)^{\nu\mp m+1}}\ln\frac{s+1}{z+1}.
\label{4.7}
\end{equation}
If in Eq.\ (\ref{4.7}) we change the \emph{name\/} of the integration
variable from $s$ to $t$, then switch from the contour
$\mathscr{C}^{\prime\prime\prime\,(-)}$ to the oppositely traversed
contour $\mathscr{C}^{\prime\prime\prime\,(+)}$ and deform the latter
into the path $\mathscr{C}^{(+)}$ (because of the structure of the
integrand, this deformation does not change the value of the integral
in question), after subsequent use of Eq.\ (\ref{3.1}), we obtain
\begin{eqnarray}
I_{\nu}^{\pm m}(z) &=& \frac{\Gamma(\nu\pm m+1)}{\Gamma(\nu+1)}
\frac{1}{2^{\nu+1}\pi\mathrm{i}}(z^{2}-1)^{\mp m/2}
\oint_{\mathscr{C}^{(+)}}\mathrm{d}t\:
\frac{(t^{2}-1)^{\nu}}{(t-z)^{\nu\mp m+1}}\ln\frac{t+1}{2}
\nonumber \\
&& -\,P_{\nu}^{\pm m}(z)\ln\frac{z+1}{2}.
\label{4.8}
\end{eqnarray}
Upon replacing the second integral on the right-hand side of Eq.\
(\ref{4.5}) by the equivalent expression given in Eq.\ (\ref{4.8}),
we arrive at the second contour-integral representation of $\partial
P_{\nu}^{\pm m}(z)/\partial\nu$:
\begin{eqnarray}
\frac{\partial P_{\nu}^{\pm m}(z)}{\partial\nu}
&=& -P_{\nu}^{\pm m}(z)\ln\frac{z+1}{2}
+[\psi(\nu\pm m+1)-\psi(\nu+1)]P_{\nu}^{\pm m}(z)
\nonumber \\
&& +\,\frac{\Gamma(\nu\pm m+1)}{\Gamma(\nu+1)}
\frac{1}{2^{\nu+1}\pi\mathrm{i}}(z^{2}-1)^{\pm m/2}
\oint_{\mathscr{C}^{(+)}}\mathrm{d}t\:
\frac{(t^{2}-1)^{\nu}}{(t-z)^{\nu\pm m+1}}\ln\frac{t+1}{2}
\nonumber \\
&& +\,\frac{\Gamma(\nu\pm m+1)}{\Gamma(\nu+1)}
\frac{1}{2^{\nu+1}\pi\mathrm{i}}(z^{2}-1)^{\mp m/2}
\oint_{\mathscr{C}^{(+)}}\mathrm{d}t\:
\frac{(t^{2}-1)^{\nu}}{(t-z)^{\nu\mp m+1}}\ln\frac{t+1}{2}.
\label{4.9}
\end{eqnarray}
It is easy to see that for $m=0$ Eq.\ (\ref{4.9}) reduces to Eq.\
(\ref{1.1}).

Next, replace in Eq.\ (\ref{4.9}) $\nu$ by $-\nu-1$. Since, by virtue
of Eq.\ (\ref{3.7}), one has
\begin{equation}
\frac{\partial P_{-\nu-1}^{\pm m}(z)}{\partial(-\nu-1)}
=-\frac{\partial P_{\nu}^{\pm m}(z)}{\partial\nu}
\label{4.10}
\end{equation}
and since it holds that
\begin{equation}
\psi(-\nu\pm m)-\psi(-\nu)=\psi(\nu\mp m+1)-\psi(\nu+1),
\label{4.11}
\end{equation}
the replacement leads to still another expression for $\partial 
P_{\nu}^{\pm m}(z)/\partial\nu$:
\begin{eqnarray}
\frac{\partial P_{\nu}^{\pm m}(z)}{\partial\nu}
&=& P_{\nu}^{\pm m}(z)\ln\frac{z+1}{2}
+[\psi(\nu+1)-\psi(\nu\mp m+1)]P_{\nu}^{\pm m}(z)
\nonumber \\
&& -\,(-)^{m}\frac{\Gamma(\nu+1)}{\Gamma(\nu\mp m+1)}
\frac{2^{\nu}}{\pi\mathrm{i}}(z^{2}-1)^{\pm m/2}
\oint_{\mathscr{C}^{(+)}}\mathrm{d}t\:
\frac{(t-z)^{\nu\mp m}}{(t^{2}-1)^{\nu+1}}\ln\frac{t+1}{2}
\nonumber \\
&& -\,(-)^{m}\frac{\Gamma(\nu+1)}{\Gamma(\nu\mp m+1)}
\frac{2^{\nu}}{\pi\mathrm{i}}(z^{2}-1)^{\mp m/2}
\oint_{\mathscr{C}^{(+)}}\mathrm{d}t\:
\frac{(t-z)^{\nu\pm m}}{(t^{2}-1)^{\nu+1}}\ln\frac{t+1}{2}.
\label{4.12}
\end{eqnarray}
Subjecting both integrals in Eq.\ (\ref{4.12}) to the variable
transformation (\ref{3.10}) results in
\begin{eqnarray}
\frac{\partial P_{\nu}^{\pm m}(z)}{\partial\nu}
&=& P_{\nu}^{\pm m}(z)\ln\frac{z+1}{2}
+[\psi(\nu+1)-\psi(\nu\mp m+1)]P_{\nu}^{\pm m}(z)
\nonumber \\
&& +\,\frac{\Gamma(\nu+1)}{\Gamma(\nu\mp m+1)}
\frac{1}{2^{\nu+1}\pi\mathrm{i}}
\left(\frac{z-1}{z+1}\right)^{\pm m/2}
\nonumber \\\
&& \quad\times\oint_{\mathscr{C}^{\prime\,(+)}}\mathrm{d}u\:
\frac{(u-1)^{\nu\mp m}(u+1)^{\nu\pm m}}{(u-z)^{\nu+1}}
\ln\frac{u+1}{z+1}
\nonumber \\
&& +\,\frac{\Gamma(\nu+1)}{\Gamma(\nu\mp m+1)}
\frac{1}{2^{\nu+1}\pi\mathrm{i}}
\left(\frac{z+1}{z-1}\right)^{\pm m/2}
\nonumber \\
&& \quad\times\oint_{\mathscr{C}^{\prime\,(+)}}\mathrm{d}u\:
\frac{(u-1)^{\nu\pm m}(u+1)^{\nu\mp m}}{(u-z)^{\nu+1}}
\ln\frac{u+1}{z+1},
\label{4.13}
\end{eqnarray}
where the path $\mathscr{C}^{\prime\,(+)}$ has been defined below
Eq.\ (\ref{3.11}). With the aid of Eqs.\ (\ref{3.11}) and
(\ref{3.16}), the above may be transformed into
\begin{eqnarray}
\frac{\partial P_{\nu}^{\pm m}(z)}{\partial\nu}
&=& -P_{\nu}^{\pm m}(z)\ln\frac{z+1}{2}
+[\psi(\nu+1)-\psi(\nu\mp m+1)]P_{\nu}^{\pm m}(z)
\nonumber \\
&& +\,\frac{\Gamma(\nu+1)}{\Gamma(\nu\mp m+1)}
\frac{1}{2^{\nu+1}\pi\mathrm{i}}
\left(\frac{z-1}{z+1}\right)^{\pm m/2}
\nonumber \\
&& \quad\times\oint_{\mathscr{C}^{\prime\,(+)}}\mathrm{d}u\:
\frac{(u-1)^{\nu\mp m}(u+1)^{\nu\pm m}}{(u-z)^{\nu+1}}
\ln\frac{u+1}{2}
\nonumber \\
&& +\,\frac{\Gamma(\nu+1)}{\Gamma(\nu\mp m+1)}
\frac{1}{2^{\nu+1}\pi\mathrm{i}}
\left(\frac{z+1}{z-1}\right)^{\pm m/2}
\nonumber \\
&& \quad\times\oint_{\mathscr{C}^{\prime\,(+)}}\mathrm{d}u\:
\frac{(u-1)^{\nu\pm m}(u+1)^{\nu\mp m}}{(u-z)^{\nu+1}}
\ln\frac{u+1}{2}.
\label{4.14}
\end{eqnarray}
%
%
\section{Formulas for $[\partial
P_{\nu}^{\pm m}(z)/\partial\nu]_{\nu=n}$}
\label{V}
\setcounter{equation}{0}
\subsection{Evaluation of $[\partial 
P_{\nu}^{m}(z)/\partial\nu]_{\nu=n}$ for $0\leqslant m\leqslant n$}
\label{V.1}
\subsubsection{Rodrigues-type formulas}
\label{V.1.1}
Let us consider Eq.\ (\ref{4.9}), with the upper signs chosen, in the
case when $\nu=n$ and $0\leqslant m\leqslant n$. We have
\begin{eqnarray}
\frac{\partial P_{\nu}^{m}(z)}{\partial\nu}\bigg|_{\nu=n}
&=& -P_{n}^{m}(z)\ln\frac{z+1}{2}
+[\psi(n+m+1)-\psi(n+1)]P_{n}^{m}(z)
\nonumber \\
&& +\,\frac{(n+m)!}{n!}\frac{1}{2^{n+1}\pi\mathrm{i}}(z^{2}-1)^{m/2}
\oint_{\mathscr{C}^{(+)}}\mathrm{d}t\:
\frac{(t^{2}-1)^{n}}{(t-z)^{n+m+1}}\ln\frac{t+1}{2}
\nonumber \\
&& +\,\frac{(n+m)!}{n!}\frac{1}{2^{n+1}\pi\mathrm{i}}(z^{2}-1)^{-m/2}
\oint_{\mathscr{C}^{(+)}}\mathrm{d}t\:
\frac{(t^{2}-1)^{n}}{(t-z)^{n-m+1}}\ln\frac{t+1}{2}
\nonumber \\
&& (0\leqslant m\leqslant n)
\label{5.1}
\end{eqnarray}
(recall that the contour $\mathscr{C}^{(+)}$ is the one defined in
Fig.\ \ref{fig1}). We see that the only singularities of the two
integrands in the domain enclosed by $\mathscr{C}^{(+)}$ are poles of
orders $n+m+1$ and $n-m+1$, respectively, located at $t=z$. Thus, on
applying the residue theorem, we find
\begin{eqnarray}
\frac{\partial P_{\nu}^{m}(z)}{\partial\nu}\bigg|_{\nu=n}
&=& -P_{n}^{m}(z)\ln\frac{z+1}{2}
+[\psi(n+m+1)-\psi(n+1)]P_{n}^{m}(z)
\nonumber \\
&& +\,\frac{1}{2^{n}n!}(z^{2}-1)^{m/2}
\frac{\mathrm{d}^{n+m}}{\mathrm{d}z^{n+m}}
\left[(z^{2}-1)^{n}\ln\frac{z+1}{2}\right]
\nonumber \\
&& +\,\frac{1}{2^{n}n!}\frac{(n+m)!}{(n-m)!}(z^{2}-1)^{-m/2}
\frac{\mathrm{d}^{n-m}}{\mathrm{d}z^{n-m}}
\left[(z^{2}-1)^{n}\ln\frac{z+1}{2}\right]
\qquad (0\leqslant m\leqslant n).
\nonumber \\
&&
\label{5.2}
\end{eqnarray}

Under the same assumptions, Eq.\ (\ref{4.14}), with the upper signs
chosen, becomes
\begin{eqnarray}
\frac{\partial P_{\nu}^{m}(z)}{\partial\nu}\bigg|_{\nu=n}
&=& -P_{n}^{m}(z)\ln\frac{z+1}{2}
+[\psi(n+1)-\psi(n-m+1)]P_{n}^{m}(z)
\nonumber \\
&& +\,\frac{n!}{(n-m)!}\frac{1}{2^{n+1}\pi\mathrm{i}}
\left(\frac{z-1}{z+1}\right)^{m/2}
\nonumber \\
&& \quad\times\oint_{\mathscr{C}^{\prime\,(+)}}\mathrm{d}u\:
\frac{(u-1)^{n-m}(u+1)^{n+m}}{(u-z)^{n+1}}\ln\frac{u+1}{2}
\nonumber \\
&& +\,\frac{n!}{(n-m)!}\frac{1}{2^{n+1}\pi\mathrm{i}}
\left(\frac{z+1}{z-1}\right)^{m/2}
\nonumber \\
&& \quad\times\oint_{\mathscr{C}^{\prime\,(+)}}\mathrm{d}u\:
\frac{(u-1)^{n+m}(u+1)^{n-m}}{(u-z)^{n+1}}\ln\frac{u+1}{2}
\qquad (0\leqslant m\leqslant n)
\label{5.3}
\end{eqnarray}
(the contour $\mathscr{C}^{\prime\,(+)}$ is the one defined below
Eq.\ (\ref{3.11})). Since in the region surrounded by
$\mathscr{C}^{\prime\,(+)}$ both integrands in Eq.\ (\ref{5.3}) have
poles of order $n+1$ located at $u=z$, by virtue of the residue
theorem we obtain
\begin{eqnarray}
\frac{\partial P_{\nu}^{m}(z)}{\partial\nu}\bigg|_{\nu=n}
&=& -P_{n}^{m}(z)\ln\frac{z+1}{2}
+[\psi(n+1)-\psi(n-m+1)]P_{n}^{m}(z)
\nonumber \\
&& +\,\frac{1}{2^{n}(n-m)!}\left(\frac{z-1}{z+1}\right)^{m/2}
\frac{\mathrm{d}^{n}}{\mathrm{d}z^{n}}
\left[(z-1)^{n-m}(z+1)^{n+m}\ln\frac{z+1}{2}\right]
\nonumber \\
&& +\,\frac{1}{2^{n}(n-m)!}\left(\frac{z+1}{z-1}\right)^{m/2}
\frac{\mathrm{d}^{n}}{\mathrm{d}z^{n}}
\left[(z-1)^{n+m}(z+1)^{n-m}\ln\frac{z+1}{2}\right]
\nonumber \\
&& (0\leqslant m\leqslant n).
\label{5.4}
\end{eqnarray}
For $m=0$, both Eqs.\ (\ref{5.2}) and (\ref{5.4}) degenerate to the
Jolliffe's formula (\ref{1.2}).
%
%
\subsubsection{Some closed-form representations}
\label{V.1.2}
The Rodrigues-type formulas (\ref{5.2}) and (\ref{5.4}) may be used
to express $[\partial P_{\nu}^{m}(z)/\partial\nu]_{\nu=n}$ with
$0\leqslant m\leqslant n$ in terms of parameter derivatives of
particular Jacobi polynomials. Using Eqs.\ (\ref{A.44}),
(\ref{A.48}), (\ref{3.36}) and (\ref{3.37}) in Eq.\ (\ref{5.2}) gives
\begin{eqnarray}
\frac{\partial P_{\nu}^{m}(z)}{\partial\nu}\bigg|_{\nu=n}
&=& P_{n}^{m}(z)\ln\frac{z+1}{2}
+[\psi(n+m+1)-\psi(n+1)]P_{n}^{m}(z)
\nonumber \\
&& +\,\frac{(n+m)!}{n!}\left(\frac{z^{2}-1}{4}\right)^{-m/2}
\frac{\partial P_{n+m}^{(-m,\beta)}(z)}
{\partial\beta}\bigg|_{\beta=-m}
\nonumber \\
&& +\,\frac{(n+m)!}{n!}\left(\frac{z^{2}-1}{4}\right)^{m/2}
\frac{\partial P_{n-m}^{(m,\beta)}(z)}
{\partial\beta}\bigg|_{\beta=m}
\qquad (0\leqslant m\leqslant n).
\label{5.5}
\end{eqnarray}
Similarly, exploiting Eqs.\ (\ref{A.29}), (\ref{A.36}), (\ref{3.32})
and (\ref{3.34}) in Eq.\ (\ref{5.4}) results in
\begin{eqnarray}
\frac{\partial P_{\nu}^{m}(z)}{\partial\nu}\bigg|_{\nu=n}
&=& P_{n}^{m}(z)\ln\frac{z+1}{2}
+[\psi(n+1)-\psi(n-m+1)]P_{n}^{m}(z)
\nonumber \\
&& +\,\frac{n!}{(n-m)!}\left(\frac{z+1}{z-1}\right)^{m/2}
\frac{\partial P_{n}^{(-m,\beta)}(z)}
{\partial\beta}\bigg|_{\beta=m}
\nonumber \\
&& +\,\frac{n!}{(n-m)!}\left(\frac{z-1}{z+1}\right)^{m/2}
\frac{\partial P_{n}^{(m,\beta)}(z)}
{\partial\beta}\bigg|_{\beta=-m}
\qquad (0\leqslant m\leqslant n).
\label{5.6}
\end{eqnarray}
The usefulness of these two formulas comes from the fact that, as it
is shown in appendix \ref{A}, it is relatively easy to obtain various
representations of the parameter derivatives of the Jacobi
polynomials entering Eqs.\ (\ref{5.5}) and (\ref{5.6}).

Using Eqs.\ (\ref{A.45}), (\ref{A.49}), (\ref{3.36}) and (\ref{3.37})
in Eq.\ (\ref{5.5}) gives
\begin{eqnarray}
\frac{\partial P_{\nu}^{m}(z)}{\partial\nu}\bigg|_{\nu=n}
&=& P_{n}^{m}(z)\ln\frac{z+1}{2}
-[\psi(n+1)+\psi(n-m+1)]P_{n}^{m}(z)
\nonumber \\
&& +\,\left(\frac{z^{2}-1}{4}\right)^{m/2}
\sum_{k=0}^{n-m}\frac{(k+n+m)!\psi(k+n+m+1)}{k!(k+m)!(n-m-k)!}
\left(\frac{z-1}{2}\right)^{k}
\nonumber \\
&& +\,\frac{(n+m)!}{(n-m)!}\left(\frac{z-1}{z+1}\right)^{m/2}
\sum_{k=0}^{n}\frac{(k+n)!\psi(k+n+1)}{k!(k+m)!(n-k)!}
\left(\frac{z-1}{2}\right)^{k}
\qquad (0\leqslant m\leqslant n).
\nonumber \\
&&
\label{5.7}
\end{eqnarray}
The same result is found if Eqs.\ (\ref{A.30}), (\ref{A.37}),
(\ref{3.32}) and (\ref{3.34}) are plugged into Eq.\ (\ref{5.6}). In
turn, inserting Eqs.\ (\ref{A.46}) and (\ref{A.50}) into Eq.\
(\ref{5.5}) and using then Eqs.\ (\ref{3.36}) and (\ref{3.37}) yields
\begin{eqnarray}
\frac{\partial P_{\nu}^{m}(z)}{\partial\nu}\bigg|_{\nu=n}
&=& P_{n}^{m}(z)\ln\frac{z+1}{2}
+[\psi(n+1)-\psi(n-m+1)]P_{n}^{m}(z)
\nonumber \\
&& -\,(-)^{n}\frac{(n+m)!}{(n-m)!}
\left(\frac{z^{2}-1}{4}\right)^{-m/2}
\nonumber \\
&& \quad \times\sum_{k=0}^{m-1}\frac{(k+n-m)!(m-k-1)!}{k!(n+m-k)!}
\left(\frac{z+1}{2}\right)^{k}
\nonumber \\
&& +\,(-)^{n+m}\left(\frac{z^{2}-1}{4}\right)^{m/2}
\sum_{k=0}^{n-m}(-)^{k}\frac{(k+n+m)!}{k!(k+m)!(n-m-k)!}
\nonumber \\
&& \quad \times[\psi(k+n+m+1)-\psi(k+m+1)]
\left(\frac{z+1}{2}\right)^{k}
\nonumber \\
&& +\,(-)^{n}\frac{(n+m)!}{(n-m)!}\left(\frac{z+1}{z-1}\right)^{m/2}
\sum_{k=0}^{n}(-)^{k}\frac{(k+n)!}{k!(k+m)!(n-k)!}
\nonumber \\
&& \quad \times[\psi(k+n+1)-\psi(k+1)]\left(\frac{z+1}{2}\right)^{k}
\qquad (0\leqslant m\leqslant n).
\label{5.8}
\end{eqnarray}
Interestingly, if Eqs.\ (\ref{A.31}), (\ref{A.39}), (\ref{3.32}) and
(\ref{3.34}) are used in Eq.\ (\ref{5.6}), one arrives at the
following representation of $[\partial
P_{\nu}^{m}(z)/\partial\nu]_{\nu=n}$ with $0\leqslant m\leqslant n$:
\begin{eqnarray}
\frac{\partial P_{\nu}^{m}(z)}{\partial\nu}\bigg|_{\nu=n} 
&=& P_{n}^{m}(z)\ln\frac{z+1}{2}
+[\psi(n+m+1)-\psi(n+1)]P_{n}^{m}(z)
\nonumber \\
&& -\,(-)^{n+m}\left(\frac{z-1}{z+1}\right)^{m/2}\sum_{k=0}^{m-1}
\frac{(k+n)!(m-k-1)!}{k!(n-k)!}\left(\frac{z+1}{2}\right)^{k}
\nonumber \\
&& +\,(-)^{n+m}\left(\frac{z^{2}-1}{4}\right)^{m/2}
\sum_{k=0}^{n-m}(-)^{k}\frac{(k+n+m)!}{k!(k+m)!(n-m-k)!}
\nonumber \\
&& \quad\times[\psi(k+n+m+1)-\psi(k+1)]
\left(\frac{z+1}{2}\right)^{k}
\nonumber \\
&& +\,(-)^{n}\frac{(n+m)!}{(n-m)!}\left(\frac{z+1}{z-1}\right)^{m/2}
\sum_{k=0}^{n}(-)^{k}\frac{(k+n)!}{k!(k+m)!(n-k)!}
\nonumber \\
&& \quad\times[\psi(k+n+1)-\psi(k+m+1)]
\left(\frac{z+1}{2}\right)^{k}
\qquad (0\leqslant m\leqslant n),
\label{5.9}
\end{eqnarray}
which does not seem to be trivially equivalent to that in Eq.\
(\ref{5.8}). Next, the formula
\begin{eqnarray}
\frac{\partial P_{\nu}^{m}(z)}{\partial\nu}\bigg|_{\nu=n} 
&=& P_{n}^{m}(z)\ln\frac{z+1}{2}
+[\psi(n+1)+\psi(n+m+1)]P_{n}^{m}(z)
\nonumber \\
&& -\,n!(n+m)!\left(\frac{z+1}{z-1}\right)^{m/2}
\left(\frac{z-1}{2}\right)^{n}
\nonumber \\
&& \quad \times\sum_{k=1}^{m}(-)^{k}
\frac{(k-1)!}{(k+n)!(k+n-m)!(m-k)!}\left(\frac{z-1}{z+1}\right)^{k}
\nonumber \\
&& -\,n!(n+m)!\left(\frac{z-1}{z+1}\right)^{m/2}
\left(\frac{z+1}{2}\right)^{n}
\nonumber \\
&& \quad \times\sum_{k=0}^{n-m}\frac{\psi(n-k+1)+\psi(n-m-k+1)}
{k!(k+m)!(n-k)!(n-m-k)!}
\left(\frac{z-1}{z+1}\right)^{k}
\qquad (0\leqslant m\leqslant n)
\nonumber \\
&&
\label{5.10}
\end{eqnarray}
is obtained if one combines Eq.\ (\ref{5.5}) with Eqs.\ (\ref{A.47}),
(\ref{A.51}), (\ref{3.36}) and (\ref{3.37}) or Eq.\ (\ref{5.6}) with
Eqs.\ (\ref{A.33}), (\ref{A.40}), (\ref{3.32}) and (\ref{3.34}). From
Eq.\ (\ref{5.10}), it is immediately found that $[\partial
P_{\nu}^{m}(z)/\partial\nu]_{\nu=n}$ with $0\leqslant m\leqslant n$
may be also written as
\begin{eqnarray}
\frac{\partial P_{\nu}^{m}(z)}{\partial\nu}\bigg|_{\nu=n} 
&=& P_{n}^{m}(z)\ln\frac{z+1}{2}
+[\psi(n+1)+\psi(n+m+1)]P_{n}^{m}(z)
\nonumber \\
&& -\,(-)^{m}n!(n+m)!\left(\frac{z-1}{z+1}\right)^{m/2}
\left(\frac{z-1}{2}\right)^{n}
\nonumber \\
&& \quad \times\sum_{k=0}^{m-1}(-)^{k}
\frac{(m-k-1)!}{k!(n-k)!(n+m-k)!}\left(\frac{z+1}{z-1}\right)^{k}
\nonumber \\
&& -\,n!(n+m)!\left(\frac{z+1}{z-1}\right)^{m/2}
\left(\frac{z-1}{2}\right)^{n}
\nonumber \\
&& \quad \times\sum_{k=0}^{n-m}\frac{\psi(k+1)+\psi(k+m+1)}
{k!(k+m)!(n-k)!(n-m-k)!}
\left(\frac{z+1}{z-1}\right)^{k}
\qquad (0\leqslant m\leqslant n).
\nonumber \\
&&
\label{5.11}
\end{eqnarray}
Finally, inserting Eqs.\ (\ref{A.35}) and (\ref{A.42}) into Eq.\
(\ref{5.6}), after subsequent use of Eqs.\ (\ref{3.32}) and
(\ref{3.34}), leads to
\begin{eqnarray}
\frac{\partial P_{\nu}^{m}(z)}{\partial\nu}\bigg|_{\nu=n}
&=& P_{n}^{m}(z)\ln\frac{z+1}{2}
+[2\psi(2n+1)-\psi(n+1)-\psi(n-m+1)]P_{n}^{m}(z)
\nonumber \\
&& +\,(-)^{n+m}\sum_{k=0}^{n-m-1}(-)^{k}
\frac{2k+2m+1}{(n-m-k)(k+n+m+1)}
\nonumber \\
&& \quad\times\left[1+\frac{k!(n+m)!}{(k+2m)!(n-m)!}\right]
P_{k+m}^{m}(z)
\nonumber \\
&& +\,(-)^{n}\frac{(n+m)!}{(n-m)!}\sum_{k=0}^{m-1}(-)^{k}
\frac{2k+1}{(n-k)(k+n+1)}P_{k}^{-m}(z)
\qquad (0\leqslant m\leqslant n).
\nonumber \\
&&
\label{5.12}
\end{eqnarray}

For $m=0$, the expressions (\ref{5.7}) to (\ref{5.12}) are seen to go
over into the representations of $[\partial
P_{\nu}(z)/\partial\nu]_{\nu=n}$ resulting from combining equation
(\ref{1.3}) with Eqs.\ (\ref{1.5}), (\ref{1.6}), (\ref{1.6}),
(\ref{1.8}), (\ref{1.9}) and (\ref{1.4}), respectively.
%
%
\subsection{Evaluation of $[\partial 
P_{\nu}^{m}(z)/\partial\nu]_{\nu=n}$ for $m>n$}
\label{V.2}
\subsubsection{Relationship between $[\partial
P_{\nu}^{m}(z)/\partial\nu]_{\nu=n}$ and $P_{n}^{-m}(z)$}
\label{V.2.1}
For $m>n$ the derivative $[\partial
P_{\nu}^{m}(z)/\partial\nu]_{\nu=n}$ may be simply related to the
function $P_{n}^{-m}(z)$. To show this, we refer to Eq.\
(\ref{3.17}), from which it follows that
\begin{equation}
P_{n}^{-m}(z)=\frac{1}{(n+m)!}
\lim_{\nu\to n}\frac{P_{\nu}^{m}(z)}{[\Gamma(\nu-m+1)]^{-1}}.
\label{5.13}
\end{equation}
Applying the l'Hospital rule and exploiting the fact that
\begin{eqnarray}
\lim_{\nu\to n}\frac{\partial}{\partial\nu}\frac{1}{\Gamma(\nu-m+1)}
&=& -\lim_{\nu\to n}\frac{\psi(\nu-m+1)}{\Gamma(\nu-m+1)}
=(-)^{n+m+1}(m-n-1)!
\qquad (m>n),
\nonumber \\
&&
\label{5.14}
\end{eqnarray}
we obtain
\begin{equation}
P_{n}^{-m}(z)=(-)^{n+m+1}\frac{1}{(n+m)!(m-n-1)!}
\frac{\partial P_{\nu}^{m}(z)}{\partial\nu}\bigg|_{\nu=n}
\qquad (m>n)
\label{5.15}
\end{equation}
and consequently
\begin{equation}
\frac{\partial P_{\nu}^{m}(z)}{\partial\nu}\bigg|_{\nu=n}
=(-)^{n+m+1}(n+m)!(m-n-1)!P_{n}^{-m}(z)
\qquad (m>n).
\label{5.16}
\end{equation}
%
%
\subsubsection{Rodrigues-type formulas}
\label{V.2.2}
The following two Rodrigues-type formulas for $[\partial
P_{\nu}^{m}(z)/\partial\nu]_{\nu=n}$ with $m>n$:
\begin{eqnarray}
\frac{\partial P_{\nu}^{m}(z)}{\partial\nu}\bigg|_{\nu=n}
&=& (-)^{n+m+1}\frac{1}{2^{n}}(m-n-1)!
\left(\frac{z+1}{z-1}\right)^{m/2}
\frac{\mathrm{d}^{n}}{\mathrm{d}z^{n}}
\left[(z-1)^{n+m}(z+1)^{n-m}\right]
\nonumber \\
&& (m>n)
\label{5.17}
\end{eqnarray}
and
\begin{eqnarray}
\frac{\partial P_{\nu}^{m}(z)}{\partial\nu}\bigg|_{\nu=n} 
&=& (-)^{n+m+1}\frac{1}{2^{n}}(m-n-1)!
\left(\frac{z-1}{z+1}\right)^{m/2}
\frac{\mathrm{d}^{n}}{\mathrm{d}z^{n}}
\left[(z-1)^{n-m}(z+1)^{n+m}\right]
\nonumber \\
&& -\,(-)^{n}2^{n+1}n!(z^{2}-1)^{m/2}
\frac{\mathrm{d}^{m-n-1}}{\mathrm{d}z^{m-n-1}}(z^{2}-1)^{-n-1}
\qquad (m>n)
\label{5.18}
\end{eqnarray}
are obtained if one combines Eq.\ (\ref{5.16}) with Eqs.\
(\ref{3.27}) and (\ref{3.30}).

A further representation of this type results from Eq.\ (\ref{4.9}).
Setting in the latter $\nu=n$, choosing the upper signs and using
Eq.\ (\ref{3.23}), we have
\begin{eqnarray}
\frac{\partial P_{\nu}^{m}(z)}{\partial\nu}\bigg|_{\nu=n}
&=& \frac{(n+m)!}{n!}\frac{1}{2^{n+1}\pi\mathrm{i}}
(z^{2}-1)^{m/2}\oint_{\mathscr{C}^{(+)}}\mathrm{d}t\:
\frac{(t^{2}-1)^{n}}{(t-z)^{n+m+1}}\ln\frac{t+1}{2}
\nonumber \\
&& +\,\frac{(n+m)!}{n!}\frac{1}{2^{n+1}\pi\mathrm{i}}
(z^{2}-1)^{-m/2}\oint_{\mathscr{C}^{(+)}}\mathrm{d}t\:
\frac{(t^{2}-1)^{n}}{(t-z)^{n-m+1}}\ln\frac{t+1}{2}
\nonumber \\
&& (m>n).
\label{5.19}
\end{eqnarray}
Since $m>n$, the integrand in the second integral in the above
equation is regular in the domain enclosed by the contour
$\mathscr{C}^{(+)}$ and thus this integral vanishes. In turn, in the
same domain the integrand in the first integral has a single pole of
order $n+m+1$ located at $t=z$, so that removing the cut (\ref{3.3})
and applying the residue theorem to this integral, we find
\begin{equation}
\frac{\partial P_{\nu}^{m}(z)}{\partial\nu}\bigg|_{\nu=n}
=\frac{1}{2^{n}n!}(z^{2}-1)^{m/2}
\frac{\mathrm{d}^{n+m}}{\mathrm{d}z^{n+m}}
\left[(z^{2}-1)^{n}\ln\frac{z+1}{2}\right]
\qquad (m>n).
\label{5.20}
\end{equation}
Since the order of differentiation is greater than the degree of the
polynomial multiplying $\ln[(z+1)/2]$, Eq.\ (\ref{5.20}) may be
transformed into
\begin{equation}
\frac{\partial P_{\nu}^{m}(z)}{\partial\nu}\bigg|_{\nu=n}
=\frac{1}{2^{n}n!}(z^{2}-1)^{m/2}
\frac{\mathrm{d}^{n+m}}{\mathrm{d}z^{n+m}}
\left[(z^{2}-1)^{n}\ln(z+1)\right]
\qquad (m>n).
\label{5.21}
\end{equation}

It is worthwhile to mention that from Eqs.\ (\ref{5.15}) and
(\ref{5.21}) one may deduce the following Rodrigues-type formula for
$P_{n}^{-m}(z)$ with $m>n$:
\begin{eqnarray}
P_{n}^{-m}(z) &=& (-)^{n+m+1}\frac{1}{2^{n}n!(n+m)!(m-n-1)!}
(z^{2}-1)^{m/2}\frac{\mathrm{d}^{n+m}}{\mathrm{d}z^{n+m}}
\left[(z^{2}-1)^{n}\ln(z+1)\right]
\nonumber \\
&& (m>n),
\label{5.22}
\end{eqnarray}
supplementing the formulas given in Eqs.\ (\ref{3.27}) and
(\ref{3.30}).
%
%
\subsection{Evaluation of $[\partial
P_{\nu}^{-m}(z)/\partial\nu]_{\nu=n}$ for $0\leqslant m\leqslant
n$} 
\label{V.3}
For $0\leqslant m\leqslant n$, the derivative $[\partial
P_{\nu}^{-m}(z)/\partial\nu]_{\nu=n}$ is most conveniently found
with the aid of the relationship (\ref{3.17}). Differentiating the
latter with respect to $\nu$ gives
\begin{equation}
\frac{\partial P_{\nu}^{-m}(z)}{\partial\nu}
=\frac{\Gamma(\nu-m+1)}{\Gamma(\nu+m+1)}
\frac{\partial P_{\nu}^{m}(z)}{\partial\nu}
+[\psi(\nu-m+1)-\psi(\nu+m+1)]P_{\nu}^{-m}(z),
\label{5.23}
\end{equation}
hence, it follows that
\begin{eqnarray}
\frac{\partial P_{\nu}^{-m}(z)}{\partial\nu}\bigg|_{\nu=n}
&=& \frac{(n-m)!}{(n+m)!}
\frac{\partial P_{\nu}^{m}(z)}{\partial\nu}\bigg|_{\nu=n}
-[\psi(n+m+1)-\psi(n-m+1)]P_{n}^{-m}(z)
\nonumber \\
&& (0\leqslant m\leqslant n)
\label{5.24}
\end{eqnarray}
(cf Eq.\ (\ref{2.10})). Equation (\ref{5.24}) may be used to obtain
various representations of $[\partial
P_{\nu}^{-m}(z)/\partial\nu]_{\nu=n}$ with $0\leqslant m\leqslant n$
directly from those derived in section \ref{V.1} for $[\partial
P_{\nu}^{m}(z)/\partial\nu]_{\nu=n}$.
%
%
\subsection{Evaluation of $[\partial 
P_{\nu}^{-m}(z)/\partial\nu]_{\nu=n}$ for $m>n$}
\label{V.4}
\subsubsection{Rodrigues-type formula}
\label{V.4.1}
Choosing in Eq.\ (\ref{4.14}) the lower signs and setting then
$\nu=n$ yields
\begin{eqnarray}
\frac{\partial P_{\nu}^{-m}(z)}{\partial\nu}\bigg|_{\nu=n}
&=& -P_{n}^{-m}(z)\ln\frac{z+1}{2}
-[\psi(n+m+1)-\psi(n+1)]P_{n}^{-m}(z)
\nonumber \\
&& +\,\frac{n!}{(n+m)!}\frac{1}{2^{n+1}\pi\mathrm{i}}
\left(\frac{z+1}{z-1}\right)^{m/2}
\nonumber \\
&& \quad\times\oint_{\mathscr{C}^{\prime\,(+)}}\mathrm{d}u\:
\frac{(u-1)^{n+m}(u+1)^{n-m}}{(u-z)^{n+1}}\ln\frac{u+1}{2}
\nonumber \\
&& +\,\frac{n!}{(n+m)!}\frac{1}{2^{n+1}\pi\mathrm{i}}
\left(\frac{z-1}{z+1}\right)^{m/2}
\nonumber \\
&& \quad\times\oint_{\mathscr{C}^{\prime\,(+)}}\mathrm{d}u\:
\frac{(u-1)^{n-m}(u+1)^{n+m}}{(u-z)^{n+1}}\ln\frac{u+1}{2}
\qquad (m>n).
\label{5.25}
\end{eqnarray}
Since the integrands in both integrals in Eq.\ (\ref{5.25}) are
single-valued in the domain enclosed by $\mathscr{C}^{\prime\,(+)}$,
in both cases the cut (\ref{3.13}) in the $u$-plane may be removed.
When this is done, it becomes possible (and, as we shall see in a
moment, also convenient) to split the second integral in Eq.\
(\ref{5.25}) into a sum of two: one over the contour
$\mathscr{C}_{z}^{\prime\,(+)}$ around the point $u=z$ in the
positive sense, with the points $u=\pm1$ left outside, and the other
over the contour $\mathscr{C}_{+1}^{\prime\,(+)}$ around the point
$u=+1$ in the positive sense, with the points $u=z$ and $u=-1$ left
outside; none of the two contours is allowed to cross the cut
(\ref{3.12}). This results in
\begin{eqnarray}
\frac{\partial P_{\nu}^{-m}(z)}{\partial\nu}\bigg|_{\nu=n}
&=& -P_{n}^{-m}(z)\ln\frac{z+1}{2}
-[\psi(n+m+1)-\psi(n+1)]P_{n}^{-m}(z)
\nonumber \\
&& +\,\frac{n!}{(n+m)!}\frac{1}{2^{n+1}\pi\mathrm{i}}
\left(\frac{z+1}{z-1}\right)^{m/2}
\nonumber \\
&& \quad\times\oint_{\mathscr{C}^{\prime\,(+)}}\mathrm{d}u\:
\frac{(u-1)^{n+m}(u+1)^{n-m}}{(u-z)^{n+1}}\ln\frac{u+1}{2}
\nonumber \\
&& +\,\frac{n!}{(n+m)!}\frac{1}{2^{n+1}\pi\mathrm{i}}
\left(\frac{z-1}{z+1}\right)^{m/2}
\nonumber \\
&& \quad\times\oint_{\mathscr{C}_{z}^{\prime\,(+)}}\mathrm{d}u\:
\frac{(u-1)^{n-m}(u+1)^{n+m}}{(u-z)^{n+1}}\ln\frac{u+1}{2}
\nonumber \\
&& +\,\frac{n!}{(n+m)!}\frac{1}{2^{n+1}\pi\mathrm{i}}
\left(\frac{z-1}{z+1}\right)^{m/2}
\nonumber \\
&& \quad\times\oint_{\mathscr{C}_{+1}^{\prime\,(+)}}\mathrm{d}u\:
\frac{(u-1)^{n-m}(u+1)^{n+m}}{(u-z)^{n+1}}\ln\frac{u+1}{2}
\qquad (m>n).
\label{5.26}
\end{eqnarray}

For a while, let us focus on the last term on the right-hand side of
the above equation, i.e., on
\begin{eqnarray}
J_{n}^{-m}(z) &=& \frac{n!}{(n+m)!}\frac{1}{2^{n+1}\pi\mathrm{i}}
\left(\frac{z-1}{z+1}\right)^{m/2}
\oint_{\mathscr{C}_{+1}^{\prime\,(+)}}\mathrm{d}u\:
\frac{(u-1)^{n-m}(u+1)^{n+m}}{(u-z)^{n+1}}\ln\frac{u+1}{2}
\nonumber \\
&& (m>n).
\label{5.27}
\end{eqnarray}
Since in the domain enclosed by $\mathscr{C}_{+1}^{\prime\,(+)}$ the
only singularity of the expression under the integral sign is the
pole of order\footnote[7]{\label{FOOT7}~At first sight, it might seem
that the pole at $u=+1$ in the integrand in Eq.\ (\ref{5.27}) is of
order $m-n$. The order is lower by one, however, due to the presence
of the factor $\ln[(u+1)/2]$.} $m-n-1$ located at $u=+1$, the
integral might be taken by evaluating a residue of the integrand at
this point. However, this method appears to be inconvenient for the
present purposes. Instead, we change the integration variable to
\begin{equation}
t=-1+2\frac{z+1}{u+1},
\label{5.28}
\end{equation}
obtaining
\begin{eqnarray}
J_{n}^{-m}(z) 
&=& (-)^{m}\frac{n!}{(n+m)!}
\frac{2^{n}}{\pi\mathrm{i}}(z^{2}-1)^{m/2}
\oint_{\mathscr{C}_{z}^{(+)}}\mathrm{d}t\:
\frac{(t-z)^{n-m}}{(t^{2}-1)^{n+1}}\ln\frac{z+1}{2}
\nonumber \\
&& -\,(-)^{m}\frac{n!}{(n+m)!}
\frac{2^{n}}{\pi\mathrm{i}}(z^{2}-1)^{m/2}
\oint_{\mathscr{C}_{z}^{(+)}}\mathrm{d}t\:
\frac{(t-z)^{n-m}}{(t^{2}-1)^{n+1}}\ln\frac{t+1}{2}
\qquad (m>n),
\nonumber \\
&&
\label{5.29}
\end{eqnarray}
where the path $\mathscr{C}_{z}^{(+)}$ runs around the point $t=z$ in
the positive sense, leaves the points $t=\pm1$ outside and does not
cross the cut (\ref{3.2}). In both integrands, the only singularity
within the domain surrounded by $\mathscr{C}_{z}^{(+)}$ is the pole
of order $m-n$ located at $t=z$, so that by the theory of residues we
obtain
\begin{eqnarray}
J_{n}^{-m}(z) 
&=& (-)^{m}\frac{2^{n+1}n!}{(n+m)!(m-n-1)!}(z^{2}-1)^{m/2}
\left[\frac{\mathrm{d}^{m-n-1}}{\mathrm{d}z^{m-n-1}}
(z^{2}-1)^{-n-1}\right]\ln\frac{z+1}{2}
\nonumber \\
&& -\,(-)^{m}\frac{2^{n+1}n!}{(n+m)!(m-n-1)!}(z^{2}-1)^{m/2}
\frac{\mathrm{d}^{m-n-1}}{\mathrm{d}z^{m-n-1}}
\left[(z^{2}-1)^{-n-1}\ln\frac{z+1}{2}\right]
\nonumber \\
&& (m>n).
\label{5.30}
\end{eqnarray}
A glance at Eq.\ (\ref{3.31}) reveals that the factor in front of
$\ln[(z+1)/2]$ in the first term on the right-hand side of Eq.\
(\ref{5.30}) equals $P_{n}^{-m}(z)-(-)^{n}P_{n}^{-m}(-z)$, i.e., we
have
\begin{eqnarray}
J_{n}^{-m}(z) 
&=& \left[P_{n}^{-m}(z)-(-)^{n}P_{n}^{-m}(-z)\right]
\ln\frac{z+1}{2}
\nonumber \\
&& -\,(-)^{m}\frac{2^{n+1}n!}{(n+m)!(m-n-1)!}(z^{2}-1)^{m/2}
\frac{\mathrm{d}^{m-n-1}}{\mathrm{d}z^{m-n-1}}
\left[(z^{2}-1)^{-n-1}\ln\frac{z+1}{2}\right]
\nonumber \\
&& (m>n).
\label{5.31}
\end{eqnarray}

We return to Eq.\ (\ref{5.26}). Evaluating the first and the second
contour integrals on its right-hand side by residues and substituting
the right-hand side of Eq.\ (\ref{5.31}) for the last term leads to
the following Rodrigues-type representation of $[\partial
P_{\nu}^{-m}(z)/\partial\nu]_{\nu=n}$ for $m>n$:
\begin{eqnarray}
\frac{\partial P_{\nu}^{-m}(z)}{\partial\nu}\bigg|_{\nu=n}
&=& (-)^{n+1}P_{n}^{-m}(-z)\ln\frac{z+1}{2}
-[\psi(n+m+1)-\psi(n+1)]P_{n}^{-m}(z)
\nonumber \\
&& +\,\frac{1}{2^{n}(n+m)!}\left(\frac{z+1}{z-1}\right)^{m/2}
\frac{\mathrm{d}^{n}}{\mathrm{d}z^{n}}
\left[(z-1)^{n+m}(z+1)^{n-m}\ln\frac{z+1}{2}\right]
\nonumber \\
&& +\,\frac{1}{2^{n}(n+m)!}\left(\frac{z-1}{z+1}\right)^{m/2}
\frac{\mathrm{d}^{n}}{\mathrm{d}z^{n}}
\left[(z-1)^{n-m}(z+1)^{n+m}\ln\frac{z+1}{2}\right]
\nonumber \\
&& -\,(-)^{m}\frac{2^{n+1}n!}{(n+m)!(m-n-1)!}(z^{2}-1)^{m/2}
\nonumber \\
&& \quad\times
\frac{\mathrm{d}^{m-n-1}}{\mathrm{d}z^{m-n-1}}
\left[(z^{2}-1)^{-n-1}\ln\frac{z+1}{2}\right]
\qquad (m>n).
\label{5.32}
\end{eqnarray}
%
%
\subsubsection{Some closed-form representations}
\label{V.4.2}
Once the Rodrigues-type representation (\ref{5.32}) is known, we may
proceed analogously as in section \ref{V.1.2}. Exploiting Eqs.\
(\ref{A.29}), (\ref{A.36}), (\ref{A.52}), (\ref{3.33}), (\ref{3.35})
and (\ref{3.38}), we transform Eq.\ (\ref{5.32}) into
\begin{eqnarray}
\frac{\partial P_{\nu}^{-m}(z)}{\partial\nu}\bigg|_{\nu=n}
&=& (-)^{n}P_{n}^{-m}(-z)\ln\frac{z+1}{2}
-[\psi(n+m+1)-\psi(n+1)]P_{n}^{-m}(z)
\nonumber \\
&& +\,\frac{n!}{(n+m)!}\left(\frac{z-1}{z+1}\right)^{m/2}
\frac{\partial P_{n}^{(m,\beta)}(z)}{\partial\beta}\bigg|_{\beta=-m}
\nonumber \\
&& +\,\frac{n!}{(n+m)!}\left(\frac{z+1}{z-1}\right)^{m/2}
\frac{\partial P_{n}^{(-m,\beta)}(z)}{\partial\beta}\bigg|_{\beta=m}
\nonumber \\
&& -\,(-)^{m}\frac{n!}{(n+m)!}\left(\frac{z^{2}-1}{4}\right)^{-m/2}
\frac{\partial P_{m-n-1}^{(-m,\beta)}(z)}
{\partial\beta}\bigg|_{\beta=-m}
\qquad (m>n).
\nonumber \\
&&
\label{5.33}
\end{eqnarray}

If in Eq.\ (\ref{5.33}) use is made of relevant representations of
the parameter derivatives of the Jacobi polynomials listed in
appendix \ref{A}, and then Eqs.\ (\ref{3.33}), (\ref{3.35}) and
(\ref{3.38}) are applied, this leads to several alternative
closed-form expressions for $[\partial
P_{\nu}^{-m}(z)/\partial\nu]_{\nu=n}$ with
$m>n$ listed below. Using the representations (\ref{A.30}),
(\ref{A.38}) and (\ref{A.53}), we obtain
\begin{eqnarray}
\frac{\partial P_{\nu}^{-m}(z)}{\partial\nu}\bigg|_{\nu=n}
&=& (-)^{n}P_{n}^{-m}(-z)\ln\frac{z+1}{2}
-(-)^{n}[\psi(n+m+1)+\psi(n+1)]P_{n}^{-m}(-z)
\nonumber \\
&& +\,\left(\frac{z-1}{z+1}\right)^{m/2}
\sum_{k=0}^{n}\frac{(k+n)!\psi(k+n+1)}{k!(k+m)!(n-k)!}
\left(\frac{z-1}{2}\right)^{k}
\nonumber \\
&& +\,\frac{(-)^{n}}{(n+m)!(m-n-1)!}
\left(\frac{z+1}{z-1}\right)^{m/2}
\nonumber \\
&& \quad \times\sum_{k=0}^{n}(-)^{k}\frac{(k+n)!(m-k-1)!\psi(k+n+1)}
{k!(n-k)!}\left(\frac{z-1}{2}\right)^{k}
\nonumber \\
&& +\,(-)^{n}\left(\frac{z^{2}-1}{4}\right)^{-m/2}
\nonumber \\
&& \quad \times\sum_{k=0}^{m-n-1}\frac{(m-k-1)!\psi(n+m-k+1)}
{k!(n+m-k)!(m-n-k-1)!}\left(\frac{z-1}{2}\right)^{k}
\qquad (m>n).
\nonumber \\
&&
\label{5.34}
\end{eqnarray}
If Eqs.\ (\ref{A.32}), (\ref{A.39}) and (\ref{A.54}) are plugged into
Eq.\ (\ref{5.33}), this results in
\begin{eqnarray}
\frac{\partial P_{\nu}^{-m}(z)}{\partial\nu}\bigg|_{\nu=n}
&=& (-)^{n}P_{n}^{-m}(-z)\ln\frac{z+1}{2}
+[\psi(m-n)-\psi(n+1)]P_{n}^{-m}(z)
\nonumber \\
&& +\,\frac{1}{(n+m)!(m-n-1)!}\left(\frac{z-1}{z+1}\right)^{m/2}
\sum_{k=0}^{n}\frac{(k+n)!(m-k-1)!}{k!(n-k)!}
\nonumber \\
&& \quad \times[\psi(k+n+1)-\psi(m-k)]\left(\frac{z+1}{2}\right)^{k}
\nonumber \\
&& -\,(-)^{n}\left(\frac{z+1}{z-1}\right)^{m/2}
\sum_{k=0}^{n}(-)^{k}\frac{(k+n)!}{k!(k+m)!(n-k)!}
\nonumber \\
&& \quad \times 
[\psi(k+m+1)-\psi(k+n+1)]\left(\frac{z+1}{2}\right)^{k}
\nonumber \\
&& -\,(-)^{m}\left(\frac{z^{2}-1}{4}\right)^{-m/2}
\sum_{k=0}^{m-n-1}(-)^{k}\frac{(m-k-1)!}{k!(n+m-k)!(m-n-k-1)!}
\nonumber \\
&& \quad \times
[\psi(n+m-k+1)-\psi(m-k)]\left(\frac{z+1}{2}\right)^{k}
\qquad (m>n).
\label{5.35}
\end{eqnarray}
In turn, use of Eqs.\ (\ref{A.34}), (\ref{A.41}) and (\ref{A.55})
gives the expression
\begin{eqnarray}
\frac{\partial P_{\nu}^{-m}(z)}{\partial\nu}\bigg|_{\nu=n}
&=& (-)^{n}P_{n}^{-m}(-z)\ln\frac{z+1}{2}
-[\psi(n+m+1)-\psi(m-n)]P_{n}^{-m}(z)
\nonumber \\
&& +\,(-)^{n}[\psi(n+m+1)+\psi(n+1)]P_{n}^{-m}(-z)
\nonumber \\
&& -\,\frac{(-)^{n}}{n!(n+m)!}\left(\frac{z+1}{z-1}\right)^{m/2}
\left(\frac{z+1}{2}\right)^{-n-1}
\nonumber \\
&& \quad \times\sum_{k=0}^{m-n-1}\frac{(k+n)!(m-k-1)!\psi(k+n+1)}
{k!(m-n-k-1)!}\left(\frac{z-1}{z+1}\right)^{k}
\nonumber \\
&& -\,\frac{n!}{(m-n-1)!}\left(\frac{z-1}{z+1}\right)^{m/2}
\left(\frac{z+1}{2}\right)^{n}
\nonumber \\
&& \quad \times\sum_{k=0}^{n}(-)^{k}\frac{(k+m-n-1)!\psi(k+m-n)}
{k!(k+m)!(n-k)!}\left(\frac{z-1}{z+1}\right)^{k}
\nonumber \\
&& -\,(-)^{n}\frac{n!}{(m-n-1)!}\left(\frac{z+1}{z-1}\right)^{m/2}
\left(\frac{z+1}{2}\right)^{n}
\nonumber \\
&& \quad \times\sum_{k=0}^{n}(-)^{k}\frac{(m-k-1)!\psi(n+m-k+1)}
{k!(n-k)!(n+m-k)!}\left(\frac{z-1}{z+1}\right)^{k}
\qquad (m>n),
\nonumber \\
&&
\label{5.36}
\end{eqnarray}
from which the counterpart representation
\begin{eqnarray}
\frac{\partial P_{\nu}^{-m}(z)}{\partial\nu}\bigg|_{\nu=n}
&=& (-)^{n}P_{n}^{-m}(-z)\ln\frac{z+1}{2}
-[\psi(n+m+1)-\psi(m-n)]P_{n}^{-m}(z)
\nonumber \\
&& +\,(-)^{n}[\psi(n+m+1)+\psi(n+1)]P_{n}^{-m}(-z)
\nonumber \\
&& -\,\frac{(-)^{n}}{n!(n+m)!}\left(\frac{z-1}{z+1}\right)^{m/2}
\left(\frac{z-1}{2}\right)^{-n-1}
\nonumber \\
&& \quad \times\sum_{k=0}^{m-n-1}\frac{(k+n)!(m-k-1)!\psi(m-k)}
{k!(m-n-k-1)!}\left(\frac{z+1}{z-1}\right)^{k}
\nonumber \\
&& -\,(-)^{n}\frac{n!}{(m-n-1)!}\left(\frac{z-1}{z+1}\right)^{m/2}
\left(\frac{z-1}{2}\right)^{n}
\nonumber \\
&& \quad \times\sum_{k=0}^{n}(-)^{k}\frac{(m-k-1)!\psi(m-k)}
{k!(n-k)!(n+m-k)!}\left(\frac{z+1}{z-1}\right)^{k}
\nonumber \\
&& -\,\frac{n!}{(m-n-1)!}\left(\frac{z+1}{z-1}\right)^{m/2}
\left(\frac{z-1}{2}\right)^{n}
\nonumber \\
&& \quad \times\sum_{k=0}^{n}(-)^{k}\frac{(k+m-n-1)!\psi(k+m+1)}
{k!(k+m)!(n-k)!}\left(\frac{z+1}{z-1}\right)^{k}
\qquad (m>n)
\nonumber \\
&&
\label{5.37}
\end{eqnarray}
follows immediately. Finally, application of Eqs.\ (\ref{A.35}),
(\ref{A.43}) and (\ref{A.56}) yields
\begin{eqnarray}
\frac{\partial P_{\nu}^{-m}(z)}{\partial\nu}\bigg|_{\nu=n}
&=& (-)^{n}P_{n}^{-m}(-z)\ln\frac{z+1}{2}
-\frac{1}{2n+1}P_{n}^{-m}(z)
\nonumber \\
&& +\,(-)^{n}[\psi(2n+2)+\psi(2n+1)-\psi(n+1)-\psi(n+m+1)]
P_{n}^{-m}(-z)
\nonumber \\
&& +\,(-)^{n}\sum_{k=0}^{n-1}(-)^{k}\frac{2k+1}{(n-k)(k+n+1)}
\nonumber \\
&& \quad \times\left[P_{k}^{-m}(z)
+(-)^{n}\frac{(k+m)!(m-k-1)!}{(n+m)!(m-n-1)!}P_{k}^{-m}(-z)\right]
\nonumber \\
&& -\,\sum_{k=1}^{m-n-1}(-)^{k}\frac{2k+2n+1}{k(k+2n+1)}
\left[P_{k+n}^{-m}(z)-(-)^{k+n}P_{k+n}^{-m}(-z)\right]
\qquad (m>n).
\nonumber \\
&&
\label{5.38}
\end{eqnarray}
%
%
\subsection{The function $[\partial 
P_{\nu}^{\pm m}(x)/\partial\nu]_{\nu=n}$}
\label{V.5}
All representations of $[\partial P_{\nu}^{\pm
m}(z)/\partial\nu]_{\nu=n}$ derived so far in this section are valid
for $z\in\mathbb{C}\setminus[-1,1]$. To obtain corresponding formulas
for $[\partial P_{\nu}^{\pm m}(x)/\partial\nu]_{\nu=n}$ with
$-1\leqslant x\leqslant1$, we may use Eq.\ (\ref{3.19}).
Differentiating the latter with respect to $\nu$ and setting then
$\nu=n$ yields
\begin{eqnarray}
\frac{\partial P_{\nu}^{\pm m}(x)}{\partial\nu}\bigg|_{\nu=n}
&=& \mathrm{e}^{\pm\mathrm{i}\pi m/2}
\frac{\partial P_{\nu}^{\pm m}(x+\mathrm{i}0)}
{\partial\nu}\bigg|_{\nu=n}
=\mathrm{e}^{\mp\mathrm{i}\pi m/2}
\frac{\partial P_{\nu}^{\pm m}(x-\mathrm{i}0)}
{\partial\nu}\bigg|_{\nu=n}
\nonumber \\
&=& \frac{1}{2}\left[\mathrm{e}^{\pm\mathrm{i}\pi m/2}
\frac{\partial P_{\nu}^{\pm m}(x+\mathrm{i}0)}
{\partial\nu}\bigg|_{\nu=n}+\mathrm{e}^{\mp\mathrm{i}\pi m/2}
\frac{\partial P_{\nu}^{\pm m}(x-\mathrm{i}0)}
{\partial\nu}\bigg|_{\nu=n}\right].
\label{5.39}
\end{eqnarray}
Particular expressions for $[\partial P_{\nu}^{\pm
m}(x)/\partial\nu]_{\nu=n}$, including the Carlson's formulas
(\ref{2.9}) and (\ref{2.11}), arise if one combines successively Eq.\
(\ref{5.39}) with the results of sections \ref{V.1} to \ref{V.4},
using Eq.\ (\ref{3.20}) and the identities
\begin{equation}
x+1\pm\mathrm{i}0=x+1,
\qquad 
x-1\pm\mathrm{i}0=\mathrm{e}^{\pm\mathrm{i}\pi}(1-x)
\qquad (-1\leqslant x\leqslant1),
\label{5.40}
\end{equation}
whenever necessary. The procedure is straightforward and therefore we
do not list here the resulting formulas.
%
%
\section{Some applications}
\label{VI}
\setcounter{equation}{0}
\subsection{Construction of the associated Legendre function of the
second kind of integer degree and order}
\label{VI.1}
In this section, we shall apply the results of section \ref{V} to
obtain several representations of the associated Legendre function of
the second kind of integer degree and order.

The following formulas:
\begin{equation}
Q_{\nu}^{m}(z)=\frac{\pi}{2}
\frac{\mathrm{e}^{\mp\mathrm{i}\pi\nu}P_{\nu}^{m}(z)
-P_{\nu}^{m}(-z)}{\sin(\pi\nu)}
\qquad (\mathrm{Im}(z)\gtrless0)
\label{6.1} 
\end{equation}
and
\begin{equation}
Q_{\nu}^{-m}(z)=\frac{\Gamma(\nu-m+1)}{\Gamma(\nu+m+1)}
Q_{\nu}^{m}(z)
\label{6.2}
\end{equation}
may serve as the definitions of the associated Legendre function of
the second kind of non-negative and negative integer order,
respectively.

In the limit $\nu\to n$, in the case of $0\leqslant m\leqslant n$,
after exploiting the l'Hospital rule, from Eq.\ (\ref{6.1}) we obtain
\begin{equation}
Q_{n}^{m}(z)=\mp\frac{1}{2}\mathrm{i}\pi P_{n}^{m}(z)
+\frac{1}{2}\frac{\partial P_{\nu}^{m}(z)}
{\partial\nu}\bigg|_{\nu=n}
-\frac{(-)^{n}}{2}\frac{\partial P_{\nu}^{m}(-z)}
{\partial\nu}\bigg|_{\nu=n}
\qquad 
(\textrm{$0\leqslant m\leqslant n$, $\mathrm{Im}(z)\gtrless0$}),
\label{6.3}
\end{equation}
while Eq.\ (\ref{6.2}) gives
\begin{equation}
Q_{n}^{-m}(z)=\frac{(n-m)!}{(n+m)!}Q_{n}^{m}(z)
\qquad (0\leqslant m\leqslant n).
\label{6.4}
\end{equation}
Inserting particular representations of $[\partial P_{\nu}^{m}(\pm
z)/\partial\nu]_{\nu=n}$ derived in section \ref{V} into the
right-hand side of Eq.\ (\ref{6.3}) and using, whenever necessary,
some of the properties of $P_{n}^{m}(z)$ featured in section
\ref{IV}, yields a variety of formulas for $Q_{n}^{\pm m}(z)$ with
$0\leqslant m\leqslant n$. If Eq.\ (\ref{5.2}) is plugged into
Eq.\ (\ref{6.3}), use is made of the property
\begin{equation}
\ln\frac{1+z}{1-z}=\ln\frac{z+1}{z-1}\pm\mathrm{i}\pi
\qquad (\textrm{Im}\,(z)\gtrless0),
\label{6.5}
\end{equation}
and the result is combined with Eq.\ (\ref{6.4}), this leads to the
following Rodrigues-type formula:
\begin{eqnarray}
Q_{n}^{\pm m}(z)
&=& -\frac{1}{2}P_{n}^{\pm m}(z)\ln\frac{z+1}{z-1}
+\frac{1}{2^{n+1}n!}(z^{2}-1)^{\pm m/2}
\frac{\mathrm{d}^{n\pm m}}{\mathrm{d}z^{n\pm m}}
\left[(z^{2}-1)^{n}\ln\frac{z+1}{z-1}\right]
\nonumber \\
&& +\,\frac{(n\pm m)!}{(n\mp m)!}
\frac{1}{2^{n+1}n!}(z^{2}-1)^{\mp m/2}
\frac{\mathrm{d}^{n\mp m}}{\mathrm{d}z^{n\mp m}}
\left[(z^{2}-1)^{n}\ln\frac{z+1}{z-1}\right]
\qquad (0\leqslant m\leqslant n).
\nonumber \\
&&
\label{6.6}
\end{eqnarray}
If Eq.\ (\ref{5.4}) is used instead of Eq.\ (\ref{5.2}), this gives
\begin{eqnarray}
Q_{n}^{\pm m}(z) 
&=& -\frac{1}{2}P_{n}^{\pm m}(z)\ln\frac{z+1}{z-1}
\nonumber \\
&& +\,\frac{1}{2^{n+1}(n\mp m)!}\left(\frac{z-1}{z+1}\right)^{m/2}
\frac{\mathrm{d}^{n}}{\mathrm{d}z^{n}}
\left[(z-1)^{n-m}(z+1)^{n+m}\ln\frac{z+1}{z-1}\right]
\nonumber \\
&& +\,\frac{1}{2^{n+1}(n\mp m)!}\left(\frac{z+1}{z-1}\right)^{m/2}
\frac{\mathrm{d}^{n}}{\mathrm{d}z^{n}}
\left[(z-1)^{n+m}(z+1)^{n-m}\ln\frac{z+1}{z-1}\right]
\nonumber \\
&& (0\leqslant m\leqslant n).
\label{6.7}
\end{eqnarray}
The latter formula has been also found by the present author, in a
different way, in \cite[section 4]{Szmy09a}. If Eq.\ (\ref{5.7}) is
employed to evaluate $[\partial P_{\nu}^{m}(z)/\partial\nu]_{\nu=n}$
and Eq.\ (\ref{5.8}) to find $[\partial
P_{\nu}^{m}(-z)/\partial\nu]_{\nu=n}$, or vice versa, from Eqs.\
(\ref{6.3}), (\ref{6.4}) and (\ref{3.22}) we obtain
\begin{equation}
Q_{n}^{\pm m}(z)=\frac{1}{2}P_{n}^{\pm m}(z)\ln\frac{z+1}{z-1}
-W_{n-1}^{\pm m}(z)
\qquad (0\leqslant m\leqslant n),
\label{6.8}
\end{equation}
with
\begin{equation}
W_{n-1}^{-m}(z)=\frac{(n-m)!}{(n+m)!}W_{n-1}^{m}(z)
\qquad (0\leqslant m\leqslant n),
\label{6.9}
\end{equation}
where
\begin{eqnarray}
W_{n-1}^{m}(z) &=& \pm\,\psi(n+1)P_{n}^{m}(z)
\nonumber \\
&& \mp\,\frac{(\pm)^{n}(\mp)^{m}}{2}\frac{(n+m)!}{(n-m)!}
\left(\frac{z^{2}-1}{4}\right)^{-m/2}
\nonumber \\
&& \quad \times\sum_{k=0}^{m-1}(\mp)^{k}
\frac{(k+n-m)!(m-k-1)!}{k!(n+m-k)!}
\left(\frac{z\mp1}{2}\right)^{k}
\nonumber \\
&& \mp\,\frac{(\pm)^{n+m}}{2}\left(\frac{z^{2}-1}{4}\right)^{m/2}
\sum_{k=0}^{n-m}(\pm)^{k}\frac{(k+n+m)!\psi(k+m+1)}{k!(k+m)!(n-m-k)!}
\left(\frac{z\mp1}{2}\right)^{k}
\nonumber \\
&& \mp\,\frac{(\pm)^{n}}{2}\frac{(n+m)!}{(n-m)!}
\left(\frac{z\mp1}{z\pm1}\right)^{m/2}
\sum_{k=0}^{n}(\pm)^{k}\frac{(k+n)!\psi(k+1)}{k!(k+m)!(n-k)!}
\left(\frac{z\mp1}{2}\right)^{k}
\nonumber \\
&& (0\leqslant m\leqslant n).
\label{6.10}
\end{eqnarray}
If Eq.\ (\ref{5.9}) is used instead of Eq.\ (\ref{5.8}), this results
in
\begin{eqnarray}
W_{n-1}^{m}(z) 
&=& \pm\,\frac{1}{2}[\psi(n+m+1)+\psi(n-m+1)]P_{n}^{m}(z)
\nonumber \\
&& \mp\,\frac{(\pm)^{n}(-)^{m}}{2}
\left(\frac{z\pm1}{z\mp1}\right)^{m/2}
\sum_{k=0}^{m-1}(\mp)^{k}\frac{(k+n)!(m-k-1)!}{k!(n-k)!}
\left(\frac{z\mp1}{2}\right)^{k}
\nonumber \\
&& \mp\,\frac{(\pm)^{n+m}}{2}\left(\frac{z^{2}-1}{4}\right)^{m/2}
\sum_{k=0}^{n-m}(\pm)^{k}\frac{(k+n+m)!\psi(k+1)}{k!(k+m)!(n-m-k)!}
\left(\frac{z\mp1}{2}\right)^{k}
\nonumber \\
&& \mp\,\frac{(\pm)^{n}}{2}\frac{(n+m)!}{(n-m)!}
\left(\frac{z\mp1}{z\pm1}\right)^{m/2}
\sum_{k=0}^{n}(\pm)^{k}\frac{(k+n)!\psi(k+m+1)}{k!(k+m)!(n-k)!}
\left(\frac{z\mp1}{2}\right)^{k}
\nonumber \\
&& (0\leqslant m\leqslant n).
\label{6.11}
\end{eqnarray}
In turn, if use is made of Eqs.\ (\ref{5.10}) and (\ref{5.11}), this
yields
\begin{eqnarray}
W_{n-1}^{m}(z) &=& \pm\,\frac{1}{2}n!(n+m)!
\left(\frac{z\pm1}{z\mp1}\right)^{m/2}
\left(\frac{z\mp1}{2}\right)^{n}
\nonumber \\
&& \quad \times\sum_{k=1}^{m}(-)^{k}
\frac{(k-1)!}{(k+n)!(k+n-m)!(m-k)!}
\left(\frac{z\mp1}{z\pm1}\right)^{k}
\nonumber \\
&& \mp\,\frac{(-)^{m}}{2}n!(n+m)!
\left(\frac{z\pm1}{z\mp1}\right)^{m/2}
\left(\frac{z\pm1}{2}\right)^{n}
\nonumber \\
&& \quad\times\sum_{k=0}^{m-1}(-)^{k}
\frac{(m-k-1)!}{k!(n-k)!(n+m-k)!}
\left(\frac{z\mp1}{z\pm1}\right)^{k}
\nonumber \\
&& \pm\,\frac{1}{2}n!(n+m)!
\left(\frac{z\mp1}{z\pm1}\right)^{m/2}
\left(\frac{z\pm1}{2}\right)^{n}
\sum_{k=0}^{n-m}\frac{1}{k!(k+m)!(n-k)!(n-m-k)!}
\nonumber \\
&& \quad \times[\psi(n-m-k+1)+\psi(n-k+1)-\psi(k+m+1)-\psi(k+1)]
\left(\frac{z\mp1}{z\pm1}\right)^{k}
\nonumber \\
&& \qquad (0\leqslant m\leqslant n).
\label{6.12}
\end{eqnarray}
Finally, application of Eq.\ (\ref{5.12}) to evaluation of both
$[\partial P_{\nu}^{m}(z)/\partial\nu]_{\nu=n}$ and $[\partial 
P_{\nu}^{m}(-z)/\partial\nu]_{\nu=n}$ leads to
\begin{eqnarray}
W_{n-1}^{m}(z) &=& \frac{1}{2}\frac{(n+m)!}{(n-m)!}
\sum_{k=0}^{m-1}(-)^{k}\frac{2k+1}{(n-k)(k+n+1)}
\left[P_{k}^{-m}(-z)-(-)^{n}P_{k}^{-m}(z)\right]
\nonumber \\
&& +\,\sum_{k=0}^{n-m-1}\frac{1-(-)^{k+n+m}}{2}
\frac{2k+2m+1}{(n-m-k)(k+n+m+1)}
\nonumber \\
&& \quad\times\left[1+\frac{k!(n+m)!}{(k+2m)!(n-m)!}\right]
P_{k+m}^{m}(z)
\qquad (0\leqslant m\leqslant n),
\label{6.13}
\end{eqnarray}
which may be more conveniently rewritten as
\begin{eqnarray}
W_{n-1}^{m}(z) &=& \frac{1}{2}\frac{(n+m)!}{(n-m)!}
\sum_{k=0}^{m-1}(-)^{k}\frac{2k+1}{(n-k)(k+n+1)}
\left[P_{k}^{-m}(-z)-(-)^{n}P_{k}^{-m}(z)\right]
\nonumber \\
&& +\,\frac{1}{2}\sum_{k=0}^{\mathrm{int}[(n-m-1)/2]}
\frac{2n-4k-1}{(n-k)(2k+1)}
\nonumber \\
&& \quad
\times\left[1+\frac{(n+m)!(n-m-2k-1)!}{(n-m)!(n+m-2k-1)!}\right]
P_{n-2k-1}^{m}(z)
\qquad (0\leqslant m\leqslant n).
\label{6.14}
\end{eqnarray}
From Eq.\ (\ref{6.9}) and either of Eqs.\ (\ref{6.10}) to
(\ref{6.14}) it is seen that the functions $W_{n-1}^{\pm m}(z)$
possess the property
\begin{equation}
W_{n-1}^{\pm m}(-z)=(-)^{n+1}W_{n-1}^{\pm m}(z)
\qquad (0\leqslant m\leqslant n).
\label{6.15}
\end{equation}

Some of the above representations of $W_{n-1}^{m}(z)$ were already
obtained, in different ways, in earlier works. In particular, the
representations in Eq.\ (\ref{6.11}) were derived by Robin
\cite[pp.\ 81, 82 and 85]{Robi58} (in this connection, cf footnote
\ref{FOOT2} on p.\ \pageref{FOOT2}), while these in Eq.\ (\ref{6.12})
may be deduced from the findings of Snow \cite[pp.\ 55 and
56]{Snow52}; an alternative method of arriving at the expressions
(\ref{6.11}) to (\ref{6.14}) has been also presented by the author in
\cite[section 4]{Szmy09a}.

We proceed to the case of $m>n$. In virtue of Eq.\ (\ref{3.23}), from
Eq.\ (\ref{6.1}) now we have
\begin{equation}
Q_{n}^{m}(z)=\frac{1}{2}\frac{\partial P_{\nu}^{m}(z)}
{\partial\nu}\bigg|_{\nu=n}
-\frac{(-)^{n}}{2}\frac{\partial P_{\nu}^{m}(-z)}
{\partial\nu}\bigg|_{\nu=n}
\qquad (m>n).
\label{6.16}
\end{equation}
Evaluating the first term on the right-hand side of Eq.\ (\ref{6.16})
with the aid of Eq.\ (\ref{5.17}) and the second one with the aid of
Eq.\ (\ref{5.18}) (or vice versa), we arrive at the Rodrigues-type
formula
\begin{equation}
Q_{n}^{m}(z)=(-)^{n+1}2^{n}n!(z^{2}-1)^{m/2}
\frac{\mathrm{d}^{m-n-1}}{\mathrm{d}z^{m-n-1}}(z^{2}-1)^{-n-1}
\qquad (m>n).
\label{6.17}
\end{equation}
Alternatively, we may use Eq.\ (\ref{5.20}) (or Eq.\ (\ref{5.21})) in
Eq.\ (\ref{6.16}). This gives
\begin{equation}
Q_{n}^{m}(z)=\frac{1}{2^{n+1}n!}(z^{2}-1)^{m/2}
\frac{\mathrm{d}^{n+m}}{\mathrm{d}z^{n+m}}
\left[(z^{2}-1)^{n}\ln\frac{1+z}{1-z}\right]
\qquad (m>n).
\label{6.18}
\end{equation}
Using Eq.\ (\ref{6.5}) and observing that in Eq.\ (\ref{6.18}) the
order of differentiation is greater than the degree of the polynomial
multiplying the logarithm, the above formula may be cast into
\begin{equation}
Q_{n}^{m}(z)=\frac{1}{2^{n+1}n!}(z^{2}-1)^{m/2}
\frac{\mathrm{d}^{n+m}}{\mathrm{d}z^{n+m}}
\left[(z^{2}-1)^{n}\ln\frac{z+1}{z-1}\right]
\qquad (m>n).
\label{6.19}
\end{equation}
Another remarkably simple expression for $Q_{n}^{m}(z)$ with $m>n$
follows if in Eq.\ (\ref{6.16}) one uses Eq.\ (\ref{5.16}); it is
\begin{equation}
Q_{n}^{m}(z)=\frac{(-)^{n+m+1}}{2}(n+m)!(m-n-1)!
\left[P_{n}^{-m}(z)-(-)^{n}P_{n}^{-m}(-z)\right]
\qquad (m>n)
\label{6.20}
\end{equation}
(cf \cite[Eq.\ (63) on p.\ 35]{Robi58}). 

The following relation:
\begin{equation}
Q_{-\nu-1}^{m}(z)=Q_{\nu}^{m}(z)-\pi\cot(\pi\nu)P_{\nu}^{m}(z)
\label{6.21}
\end{equation}
may be easily derived from Eq.\ (\ref{6.1}). From it, by virtue of
Eq.\ (\ref{3.23}), one finds
\begin{equation}
Q_{-n-1}^{m}(z)=Q_{n}^{m}(z)
-\frac{\partial P_{\nu}^{m}(z)}{\partial\nu}\bigg|_{\nu=n}
\qquad (m>n).
\label{6.22}
\end{equation}
It is thus seen that representations of $Q_{-n-1}^{m}(z)$ with $m>n$
may be straightforwardly deduced from those of $Q_{n}^{m}(z)$, with
the use of the findings of section \ref{V.2}. In this way, one
arrives at
\begin{equation}
Q_{-n-1}^{m}(z)=\frac{(-)^{n+m}}{2}(n+m)!(m-n-1)!
\left[P_{n}^{-m}(z)+(-)^{n}P_{n}^{-m}(-z)\right]
\qquad (m>n)
\label{6.23}
\end{equation}
and
\begin{equation}
Q_{-n-1}^{m}(z)=-\frac{1}{2^{n+1}n!}(z^{2}-1)^{m/2}
\frac{\mathrm{d}^{n+m}}{\mathrm{d}z^{n+m}}
\left[(z^{2}-1)^{n}\ln(z^{2}-1)\right]
\qquad (m>n).
\label{6.24}
\end{equation}
Other Rodrigues-type formulas follow if one couples Eq.\ (\ref{6.23})
with Eqs.\ (\ref{3.27}) and (\ref{3.30}).

Concluding, we observe that on the cut $-1\leqslant x\leqslant1$
counterpart expressions for the associated Legendre function of the
second kind of integer order may be obtained from the results of this
section by combining them with the defining formula
\begin{equation}
Q_{\nu}^{\pm m}(x)=\frac{(-)^{m}}{2}
\left[\mathrm{e}^{\mp\mathrm{i}\pi m/2}
Q_{\nu}^{\pm m}(x+\mathrm{i}0)
+\mathrm{e}^{\pm\mathrm{i}\pi m/2}
Q_{\nu}^{\pm m}(x-\mathrm{i}0)\right].
\label{6.25}
\end{equation}
%
%
\subsection{Evaluation of 
$[\partial^{2}P_{\nu}^{m}(z)/\partial\nu^{2}]_{\nu=n}$ for $m>n$}
\label{VI.2}
In this section, we shall show that if $m>n$, then the knowledge of
$[\partial P_{\nu}^{-m}(z)/\partial\nu]_{\nu=n}$ allows one to
evaluate $[\partial^{2}P_{\nu}^{m}(z)/\partial\nu^{2}]_{\nu=n}$.

To begin, we observe that from the easily provable (cf Eqs.\
(\ref{2.3}) and (\ref{2.4})) identity
\begin{equation}
\psi(\zeta)=\psi(1-\zeta)-\cos(\pi\zeta)\Gamma(\zeta)\Gamma(1-\zeta)
\label{6.26}
\end{equation}
it follows that
\begin{equation}
\psi(\nu-m+1)=\psi(m-\nu)
-\cos[\pi(\nu-m+1)]\Gamma(m-\nu)\Gamma(\nu-m+1).
\label{6.27}
\end{equation}
With this, Eq.\ (\ref{5.23}) may be rewritten as
\begin{eqnarray}
\frac{\partial P_{\nu}^{-m}(z)}{\partial\nu}
&=& [\psi(m-\nu)-\psi(\nu+m+1)]P_{\nu}^{-m}(z)
+\frac{\Gamma(\nu-m+1)}{\Gamma(\nu+m+1)}
\nonumber \\
&& \times\left\{\frac{\partial P_{\nu}^{m}(z)}{\partial\nu}
-\cos[\pi(\nu-m+1)]\Gamma(m-\nu)\Gamma(\nu+m+1)P_{\nu}^{-m}(z)
\right\}.
\label{6.28}
\end{eqnarray}
In the limit $\nu\to n$ (with $m>n$), Eq.\ (\ref{6.28}) becomes
\begin{eqnarray}
\frac{\partial P_{\nu}^{-m}(z)}{\partial\nu}\bigg|_{\nu=n}
&=& -\,[\psi(n+m+1)-\psi(m-n)]P_{n}^{-m}(z)
\nonumber \\
&& +\,\frac{1}{(n+m)!}\lim_{\nu\to n}\Gamma(\nu-m+1)
\bigg\{\frac{\partial P_{\nu}^{m}(z)}{\partial\nu}
\nonumber \\
&& \quad -\,\cos[\pi(\nu-m+1)]\Gamma(m-\nu)\Gamma(\nu+m+1)
P_{\nu}^{-m}(z)\bigg\}
\qquad (m>n).
\nonumber \\
&&
\label{6.29}
\end{eqnarray}
The limit which remains to be evaluated on the right-hand side of
Eq.\ (\ref{6.29}) may be taken with the aid of the l'Hospital rule.
This gives
\begin{eqnarray}
\frac{\partial P_{\nu}^{-m}(z)}{\partial\nu}\bigg|_{\nu=n}
&=& -2[\psi(n+m+1)-\psi(m-n)]P_{n}^{-m}(z)
-\frac{\partial P_{\nu}^{-m}(z)}{\partial\nu}\bigg|_{\nu=n}
\nonumber \\
&& -\,(-)^{n+m}\frac{1}{(n+m)!(m-n-1)!}
\frac{\partial^{2}P_{\nu}^{m}(z)}{\partial\nu^{2}}\bigg|_{\nu=n}
\qquad (m>n).
\label{6.30}
\end{eqnarray}
Solving Eq.\ (\ref{6.30}) for
$[\partial^{2}P_{\nu}^{m}(z)/\partial\nu^{2}]_{\nu=n}$ results in
the sought relationship:
\begin{eqnarray}
\frac{\partial^{2}P_{\nu}^{m}(z)}{\partial\nu^{2}}\bigg|_{\nu=n}
&=& (-)^{n+m+1}2(n+m)!(m-n-1)!
\nonumber \\
&& \times\left\{[\psi(n+m+1)-\psi(m-n)]P_{n}^{-m}(z)
+\frac{\partial P_{\nu}^{-m}(z)}
{\partial\nu}\bigg|_{\nu=n}\right\}
\qquad (m>n).
\nonumber \\
&&
\label{6.31}
\end{eqnarray}
Various explicit representations of 
$[\partial^{2}P_{\nu}^{m}(z)/\partial\nu^{2}]_{\nu=n}$ with $m>n$ may
be obtained from this formula with the aid of the results of section
\ref{V.4}.

A counterpart expression for
$[\partial^{2}P_{\nu}^{m}(x)/\partial\nu^{2}]_{\nu=n}$ with $m>n$,
in terms of $[\partial P_{\nu}^{-m}(x)/\partial\nu]_{\nu=n}$ and 
$P_{n}^{-m}(x)$, follows if one combines Eq.\ (\ref{6.31}) with Eqs.\ 
(\ref{3.19}) and (\ref{5.39}).
%
%
\subsection{Evaluation of $[\partial 
Q_{\nu}^{m}(z)/\partial\nu]_{\nu=n}$ and $[\partial 
Q_{\nu}^{m}(z)/\partial\nu]_{\nu=-n-1}$ for $m>n$}
\label{VI.3}
Finally, below we shall show that for $m>n$ it is possible to relate
the derivatives
$[\partial Q_{\nu}^{m}(z)/\partial\nu]_{\nu=n}$ and $[\partial
Q_{\nu}^{m}(z)/\partial\nu]_{\nu=-n-1}$ to the derivatives
$[\partial P_{\nu}^{-m}(\pm z)/\partial\nu]_{\nu=n}$.

Differentiating Eq.\ (\ref{6.1}) with respect to $\nu$ gives
\begin{eqnarray}
\frac{\partial Q_{\nu}^{m}(z)}{\partial\nu}
&=& \frac{\pi}{\sin(\pi\nu)}\Bigg[-\cos(\pi\nu)Q_{\nu}^{m}(z)
\mp\frac{1}{2}\mathrm{i}\pi\mathrm{e}^{\mp\mathrm{i}\pi\nu}
P_{\nu}^{m}(z)
\nonumber \\
&& +\,\frac{1}{2}\mathrm{e}^{\mp\mathrm{i}\pi\nu}
\frac{\partial P_{\nu}^{m}(z)}{\partial\nu}
-\frac{1}{2}\frac{\partial P_{\nu}^{m}(-z)}{\partial\nu}\Bigg]
\qquad (\mathrm{Im}(z)\gtrless0).
\label{6.32}
\end{eqnarray}
In the limit $\nu\to n$, after making use of the l'Hospital rule,
from Eq.\ (\ref{6.32}) we have
\begin{eqnarray}
\frac{\partial Q_{\nu}^{m}(z)}{\partial\nu}\bigg|_{\nu=n}
&=& -\frac{\partial Q_{\nu}^{m}(z)}{\partial\nu}\bigg|_{\nu=n}
-\frac{1}{2}\pi^{2}P_{n}^{m}(z)
\mp\mathrm{i}\pi\frac{\partial P_{\nu}^{m}(z)}
{\partial\nu}\bigg|_{\nu=n}
\nonumber \\
&& +\,\frac{1}{2}\frac{\partial^{2}P_{\nu}^{m}(z)}
{\partial\nu^{2}}\bigg|_{\nu=n}
-\frac{(-)^{n}}{2}\frac{\partial^{2}P_{\nu}^{m}(-z)}
{\partial\nu^{2}}\bigg|_{\nu=n}
\qquad (\mathrm{Im}(z)\gtrless0).
\label{6.33}
\end{eqnarray}
Solving Eq.\ (\ref{6.33}) for $[\partial
Q_{\nu}^{m}(z)/\partial\nu]_{\nu=n}$ gives
\begin{eqnarray}
\frac{\partial Q_{\nu}^{m}(z)}{\partial\nu}\bigg|_{\nu=n}
&=& -\,\frac{1}{4}\pi^{2}P_{n}^{m}(z)
\mp\frac{1}{2}\mathrm{i}\pi\frac{\partial P_{\nu}^{m}(z)}
{\partial\nu}\bigg|_{\nu=n}
+\frac{1}{4}\frac{\partial^{2}P_{\nu}^{m}(z)}
{\partial\nu^{2}}\bigg|_{\nu=n}
\nonumber \\
&& -\,\frac{(-)^{n}}{4}\frac{\partial^{2}P_{\nu}^{m}(-z)}
{\partial\nu^{2}}\bigg|_{\nu=n}
\qquad (\mathrm{Im}(z)\gtrless0).
\label{6.34}
\end{eqnarray}
So far, $n$ has been an arbitrary non-negative integer. Imposing in
Eq.\ (\ref{6.34}) the restriction $m>n$ and using Eq.\ (\ref{3.23}),
we obtain
\begin{eqnarray}
\frac{\partial Q_{\nu}^{m}(z)}{\partial\nu}\bigg|_{\nu=n}
&=& \mp\,\frac{1}{2}\mathrm{i}\pi\frac{\partial P_{\nu}^{m}(z)}
{\partial\nu}\bigg|_{\nu=n}
+\frac{1}{4}\frac{\partial^{2}P_{\nu}^{m}(z)}
{\partial\nu^{2}}\bigg|_{\nu=n}
-\frac{(-)^{n}}{4}\frac{\partial^{2}P_{\nu}^{m}(-z)}
{\partial\nu^{2}}\bigg|_{\nu=n}
\nonumber \\
&& (\textrm{$m>n$, $\mathrm{Im}(z)\gtrless0$}).
\label{6.35}
\end{eqnarray}
To eliminate the second derivatives from the right-hand side of Eq.\
(\ref{6.35}), we may exploit Eq.\ (\ref{6.31}); in this way, using
additionally Eqs.\ (\ref{5.16}) and (\ref{6.20}), we find
\begin{eqnarray}
\frac{\partial Q_{\nu}^{m}(z)}{\partial\nu}\bigg|_{\nu=n}
&=& [\psi(n+m+1)-\psi(m-n)]Q_{n}^{m}(z)
\nonumber \\
&& +\,\frac{(-)^{n+m}}{2}(n+m)!(m-n-1)!
\nonumber \\
&& \quad \times\left[\pm\mathrm{i}\pi P_{n}^{-m}(z)
-\frac{\partial P_{\nu}^{-m}(z)}{\partial\nu}\bigg|_{\nu=n}
+(-)^{n}\frac{\partial P_{\nu}^{-m}(-z)}{\partial\nu}\bigg|_{\nu=n}
\right]
\nonumber \\
&& (\textrm{$m>n$, $\mathrm{Im}(z)\gtrless0$}).
\label{6.36}
\end{eqnarray}
Proceeding analogously, with the aid of Eq.\ (\ref{6.23}) and the
relationship
\begin{equation}
\frac{\partial^{2}P_{\nu}^{m}(z)}{\partial\nu^{2}}\bigg|_{\nu=-n-1}
=\frac{\partial^{2}P_{\nu}^{m}(z)}{\partial\nu^{2}}\bigg|_{\nu=n},
\label{6.37}
\end{equation}
resulting from Eq.\ (\ref{3.7}), one obtains
\begin{eqnarray}
\frac{\partial Q_{\nu}^{m}(z)}{\partial\nu}\bigg|_{\nu=-n-1}
&=& -\,[\psi(n+m+1)-\psi(m-n)]Q_{-n-1}^{m}(z)
\nonumber \\
&& -\,\frac{(-)^{n+m}}{2}(n+m)!(m-n-1)!
\nonumber \\
&& \quad \times\left[\pm\mathrm{i}\pi P_{n}^{-m}(z)
+\frac{\partial P_{\nu}^{-m}(z)}{\partial\nu}\bigg|_{\nu=n}
+(-)^{n}\frac{\partial P_{\nu}^{-m}(-z)}{\partial\nu}\bigg|_{\nu=n}
\right]
\nonumber \\
&& (\textrm{$m>n$, $\mathrm{Im}(z)\gtrless0$}).
\label{6.38}
\end{eqnarray}
Equations (\ref{6.36}) and (\ref{6.38}) may be combined with the
formulas found in section \ref{V.4} to yield several explicit
representations of $[\partial Q_{\nu}^{m}(z)/\partial\nu]_{\nu=n}$
and $[\partial Q_{\nu}^{m}(z)/\partial\nu]_{\nu=-n-1}$ with $m>n$.

To derive counterpart expressions on the cut $-1\leqslant
x\leqslant1$, one should use the results (\ref{6.36}) and
(\ref{6.38}) in conjunction with Eqs.\ (\ref{6.25}), (\ref{3.19}) and
(\ref{5.39}).
\section*{Acknowledgments}
The author wishes to thank an anonymous referee to \cite{Szmy06b},
whose suggestion to use the contour-integration technique to evaluate
the derivative $[\partial P_{\nu}(z)/\partial\nu]_{\nu=n}$ inspired
the present work.
%
%
\appendix
\section{Appendix: Some relevant properties of the Jacobi polynomials}
\label{A}
\setcounter{equation}{0}
The Jacobi polynomials \cite{Szeg39} may be defined through the
Rodrigues-type formula
\begin{equation}
P_{n}^{(\alpha,\beta)}(z)=\frac{1}{2^{n}n!}
(z-1)^{-\alpha}(z+1)^{-\beta}
\frac{\mathrm{d}^{n}}{\mathrm{d}z^{n}}
\left[(z-1)^{n+\alpha}(z+1)^{n+\beta}\right]
\qquad (\alpha,\beta\in\mathbb{C}).
\label{A.1}
\end{equation}
If
\begin{equation}
\frac{\Gamma(2n+\alpha+\beta+1)}{\Gamma(n+\alpha+\beta+1)}\neq0,
\label{A.2}
\end{equation}
then $P_{n}^{(\alpha,\beta)}(z)$ is a polynomial in $z$ of degree
$n$. From Eq.\ (\ref{A.1}) it is seen that
\begin{equation}
P_{n}^{(\beta,\alpha)}(-z)=(-)^{n}P_{n}^{(\alpha,\beta)}(z).
\label{A.3}
\end{equation}

The following explicit representations of $P_{n}^{(\alpha,\beta)}(z)$
have proved to be useful in the context of the present paper:
\begin{equation}
P_{n}^{(\alpha,\beta)}(z)
=\frac{\Gamma(n+\alpha+1)}{\Gamma(n+\alpha+\beta+1)}
\sum_{k=0}^{n}\frac{\Gamma(k+n+\alpha+\beta+1)}
{k!(n-k)!\Gamma(k+\alpha+1)}\left(\frac{z-1}{2}\right)^{k},
\label{A.4}
\end{equation}
\begin{eqnarray}
P_{n}^{(\alpha,\beta)}(z)
&=& \frac{(-)^{n}}{\Gamma(-n-\alpha)\Gamma(n+\alpha+\beta+1)}
\nonumber \\
&& \times \sum_{k=0}^{n}(-)^{k}
\frac{\Gamma(-k-\alpha)\Gamma(k+n+\alpha+\beta+1)}{k!(n-k)!}
\left(\frac{z-1}{2}\right)^{k},
\label{A.5}
\end{eqnarray}
\begin{equation}
P_{n}^{(\alpha,\beta)}(z)
=(-)^{n}\frac{\Gamma(-n-\alpha-\beta)}{\Gamma(-n-\alpha)}
\sum_{k=0}^{n}\frac{\Gamma(-k-\alpha)}
{k!(n-k)!\Gamma(-k-n-\alpha-\beta)}\left(\frac{z-1}{2}\right)^{k},
\label{A.6}
\end{equation}
\begin{equation}
P_{n}^{(\alpha,\beta)}(z)
=(-)^{n}\frac{\Gamma(n+\beta+1)}{\Gamma(n+\alpha+\beta+1)}
\sum_{k=0}^{n}(-)^{k}\frac{\Gamma(k+n+\alpha+\beta+1)}
{k!(n-k)!\Gamma(k+\beta+1)}\left(\frac{z+1}{2}\right)^{k},
\label{A.7}
\end{equation}
\begin{eqnarray}
P_{n}^{(\alpha,\beta)}(z)
&=& \frac{1}{\Gamma(-n-\beta)\Gamma(n+\alpha+\beta+1)}
\nonumber \\
&& \times\sum_{k=0}^{n}
\frac{\Gamma(-k-\beta)\Gamma(k+n+\alpha+\beta+1)}
{k!(n-k)!}\left(\frac{z+1}{2}\right)^{k},
\label{A.8}
\end{eqnarray}
\begin{equation}
P_{n}^{(\alpha,\beta)}(z)
=\frac{\Gamma(-n-\alpha-\beta)}{\Gamma(-n-\beta)}
\sum_{k=0}^{n}(-)^{k}\frac{\Gamma(-k-\beta)}
{k!(n-k)!\Gamma(-k-n-\alpha-\beta)}\left(\frac{z+1}{2}\right)^{k},
\label{A.9}
\end{equation}
\begin{eqnarray}
P_{n}^{(\alpha,\beta)}(z) &=& \Gamma(n+\alpha+1)\Gamma(n+\beta+1)
\left(\frac{z+1}{2}\right)^{n}
\nonumber \\
&& \times\sum_{k=0}^{n}
\frac{1}{k!(n-k)!\Gamma(k+\alpha+1)\Gamma(n+\beta-k+1)}
\left(\frac{z-1}{z+1}\right)^{k},
\label{A.10}
\end{eqnarray}
\begin{equation}
P_{n}^{(\alpha,\beta)}(z) 
=\frac{\Gamma(n+\alpha+1)}{\Gamma(-n-\beta)}
\left(\frac{z+1}{2}\right)^{n}\sum_{k=0}^{n}(-)^{k}
\frac{\Gamma(k-n-\beta)}{k!(n-k)!\Gamma(k+\alpha+1)}
\left(\frac{z-1}{z+1}\right)^{k},
\label{A.11}
\end{equation}
\begin{eqnarray}
P_{n}^{(\alpha,\beta)}(z) 
&=& (-)^{n}\frac{\Gamma(n+\beta+1)}{\Gamma(-n-\alpha)}
\left(\frac{z+1}{2}\right)^{n}
\nonumber \\ 
&& \times\sum_{k=0}^{n}(-)^{k}
\frac{\Gamma(-k-\alpha)}{k!(n-k)!\Gamma(n+\beta-k+1)}
\left(\frac{z-1}{z+1}\right)^{k}
\label{A.12}
\end{eqnarray}
and
\begin{eqnarray}
P_{n}^{(\alpha,\beta)}(z)
&=& \frac{(-)^{n}}{\Gamma(-n-\alpha)\Gamma(-n-\beta)}
\left(\frac{z+1}{2}\right)^{n}\sum_{k=0}^{n}
\frac{\Gamma(-k-\alpha)\Gamma(k-n-\beta)}{k!(n-k)!}
\left(\frac{z-1}{z+1}\right)^{k}.
\nonumber \\
&&
\label{A.13}
\end{eqnarray}

Differentiation of Eq.\ (\ref{A.1}) with respect to $\beta$ gives
\begin{eqnarray}
\frac{\partial P_{n}^{(\alpha,\beta)}(z)}{\partial\beta}
&=& -\,P_{n}^{(\alpha,\beta)}(z)\ln(z+1)
\nonumber \\
&& +\,\frac{1}{2^{n}n!}(z-1)^{-\alpha}(z+1)^{-\beta}
\frac{\mathrm{d}^{n}}{\mathrm{d}z^{n}}
\left[(z-1)^{n+\alpha}(z+1)^{n+\beta}\ln(z+1)\right].
\label{A.14}
\end{eqnarray}
By virtue of Eq.\ (\ref{A.1}), the above result may be rewritten in
the form
\begin{eqnarray}
\frac{\partial P_{n}^{(\alpha,\beta)}(z)}{\partial\beta}
&=& -\,P_{n}^{(\alpha,\beta)}(z)\ln\frac{z+1}{2}
\nonumber \\
&& +\,\frac{1}{2^{n}n!}(z-1)^{-\alpha}(z+1)^{-\beta}
\frac{\mathrm{d}^{n}}{\mathrm{d}z^{n}}
\left[(z-1)^{n+\alpha}(z+1)^{n+\beta}\ln\frac{z+1}{2}\right],
\label{A.15}
\end{eqnarray}
which is more suitable for the purposes of the present paper. Next,
differentiation of Eqs.\ (\ref{A.4}) to (\ref{A.13}) with respect to
$\beta$ yields the following formulas:
\begin{eqnarray}
\frac{\partial P_{n}^{(\alpha,\beta)}(z)}{\partial\beta}
&=& -\,\psi(n+\alpha+\beta+1)P_{n}^{(\alpha,\beta)}(z)
\nonumber \\
&& +\,\frac{\Gamma(n+\alpha+1)}{\Gamma(n+\alpha+\beta+1)}
\sum_{k=0}^{n}\frac{\Gamma(k+n+\alpha+\beta+1)
\psi(k+n+\alpha+\beta+1)}{k!(n-k)!\Gamma(k+\alpha+1)}
\left(\frac{z-1}{2}\right)^{k},
\nonumber \\
&&
\label{A.16}
\end{eqnarray}
\begin{eqnarray}
\frac{\partial P_{n}^{(\alpha,\beta)}(z)}{\partial\beta} 
&=& -\,\psi(n+\alpha+\beta+1)P_{n}^{(\alpha,\beta)}(z)
\nonumber \\
&& +\,\frac{(-)^{n}}{\Gamma(-n-\alpha)\Gamma(n+\alpha+\beta+1)}
\nonumber \\
&& \quad \times\sum_{k=0}^{n}(-)^{k}
\frac{\Gamma(-k-\alpha)\Gamma(k+n+\alpha+\beta+1)
\psi(k+n+\alpha+\beta+1)}{k!(n-k)!}\left(\frac{z-1}{2}\right)^{k},
\nonumber \\
&&
\label{A.17}
\end{eqnarray}
\begin{eqnarray}
\frac{\partial P_{n}^{(\alpha,\beta)}(z)}{\partial\beta} 
&=& -\,\psi(-n-\alpha-\beta)P_{n}^{(\alpha,\beta)}(z)
\nonumber \\
&& +\,(-)^{n}\frac{\Gamma(-n-\alpha-\beta)}{\Gamma(-n-\alpha)}
\sum_{k=0}^{n}\frac{\Gamma(-k-\alpha)\psi(-k-n-\alpha-\beta)}
{k!(n-k)!\Gamma(-k-n-\alpha-\beta)}\left(\frac{z-1}{2}\right)^{k},
\nonumber \\
&&
\label{A.18}
\end{eqnarray}
\begin{eqnarray}
\frac{\partial P_{n}^{(\alpha,\beta)}(z)}{\partial\beta}
&=& [\psi(n+\beta+1)-\psi(n+\alpha+\beta+1)]
P_{n}^{(\alpha,\beta)}(z)
\nonumber \\
&& +\,(-)^{n}\frac{\Gamma(n+\beta+1)}{\Gamma(n+\alpha+\beta+1)}
\sum_{k=0}^{n}(-)^{k}\frac{\Gamma(k+n+\alpha+\beta+1)}
{k!(n-k)!\Gamma(k+\beta+1)}
\nonumber \\
&& \quad\times[\psi(k+n+\alpha+\beta+1)-\psi(k+\beta+1)]
\left(\frac{z+1}{2}\right)^{k},
\label{A.19}
\end{eqnarray}
\begin{eqnarray}
\frac{\partial P_{n}^{(\alpha,\beta)}(z)}{\partial\beta}
&=& [\psi(-n-\beta)-\psi(n+\alpha+\beta+1)]
P_{n}^{(\alpha,\beta)}(z)
\nonumber \\
&& +\,\frac{1}{\Gamma(-n-\beta)\Gamma(n+\alpha+\beta+1)}
\sum_{k=0}^{n}\frac{\Gamma(-k-\beta)\Gamma(k+n+\alpha+\beta+1)}
{k!(n-k)!}
\nonumber \\
&& \quad \times[\psi(k+n+\alpha+\beta+1)-\psi(-k-\beta)]
\left(\frac{z+1}{2}\right)^{k},
\label{A.20}
\end{eqnarray}
\begin{eqnarray}
\frac{\partial P_{n}^{(\alpha,\beta)}(z)}{\partial\beta} 
&=& [\psi(-n-\beta)-\psi(-n-\alpha-\beta)]P_{n}^{(\alpha,\beta)}(z)
\nonumber \\
&& +\,\frac{\Gamma(-n-\alpha-\beta)}{\Gamma(-n-\beta)}
\sum_{k=0}^{n}(-)^{k}\frac{\Gamma(-k-\beta)}
{k!(n-k)!\Gamma(-k-n-\alpha-\beta)}
\nonumber \\
&& \quad \times[\psi(-k-n-\alpha-\beta)-\psi(-k-\beta)]
\left(\frac{z+1}{2}\right)^{k},
\label{A.21}
\end{eqnarray}
\begin{eqnarray}
\frac{\partial P_{n}^{(\alpha,\beta)}(z)}{\partial\beta}
&=& \psi(n+\beta+1)P_{n}^{(\alpha,\beta)}(z)
\nonumber \\
&& -\,\Gamma(n+\alpha+1)\Gamma(n+\beta+1)
\left(\frac{z+1}{2}\right)^{n}
\nonumber \\
&& \quad \times\sum_{k=0}^{n}\frac{\psi(n+\beta-k+1)}
{k!(n-k)!\Gamma(k+\alpha+1)\Gamma(n+\beta-k+1)}
\left(\frac{z-1}{z+1}\right)^{k},
\label{A.22}
\end{eqnarray}
\begin{eqnarray}
\frac{\partial P_{n}^{(\alpha,\beta)}(z)}{\partial\beta}
&=& \psi(-n-\beta)P_{n}^{(\alpha,\beta)}(z)
\nonumber \\
&& -\,\frac{\Gamma(n+\alpha+1)}{\Gamma(-n-\beta)}
\left(\frac{z+1}{2}\right)^{n}
\nonumber \\
&& \quad \times\sum_{k=0}^{n}(-)^{k}
\frac{\Gamma(k-n-\beta)\psi(k-n-\beta)}
{k!(n-k)!\Gamma(k+\alpha+1)}\left(\frac{z-1}{z+1}\right)^{k},
\label{A.23}
\end{eqnarray}
\begin{eqnarray}
\frac{\partial P_{n}^{(\alpha,\beta)}(z)}{\partial\beta}
&=& \psi(n+\beta+1)P_{n}^{(\alpha,\beta)}(z)
\nonumber \\
&& -\,(-)^{n}\frac{\Gamma(n+\beta+1)}{\Gamma(-n-\alpha)}
\left(\frac{z+1}{2}\right)^{n}
\nonumber \\
&& \quad \times\sum_{k=0}^{n}(-)^{k}
\frac{\Gamma(-k-\alpha)\psi(n+\beta-k+1)}
{k!(n-k)!\Gamma(n+\beta-k+1)}\left(\frac{z-1}{z+1}\right)^{k},
\label{A.24}
\end{eqnarray}
\begin{eqnarray}
\frac{\partial P_{n}^{(\alpha,\beta)}(z)}{\partial\beta}
&=& \psi(-n-\beta)P_{n}^{(\alpha,\beta)}(z)
\nonumber \\
&& -\,\frac{(-)^{n}}{\Gamma(-n-\alpha)\Gamma(-n-\beta)}
\left(\frac{z+1}{2}\right)^{n}
\nonumber \\
&& \quad \times\sum_{k=0}^{n}
\frac{\Gamma(-k-\alpha)\Gamma(k-n-\beta)\psi(k-n-\beta)}{k!(n-k)!}
\left(\frac{z-1}{z+1}\right)^{k}.
\label{A.25}
\end{eqnarray}

Consider next the relationship \cite[equation (1.23.2.2)]{Bryc08}
(for two different proofs, see \cite{Froh94,Szmy09b})
\begin{eqnarray}
\frac{\partial P_{n}^{(\alpha,\beta)}(z)}{\partial\beta}
&=& [\psi(2n+\alpha+\beta+1)-\psi(n+\alpha+\beta+1)]
P_{n}^{(\alpha,\beta)}(z)
\nonumber \\
&& +\,(-)^{n}\frac{\Gamma(n+\alpha+1)}{\Gamma(n+\alpha+\beta+1)}
\sum_{k=0}^{n-1}(-)^{k}
\frac{2k+\alpha+\beta+1}{(n-k)(k+n+\alpha+\beta+1)}
\nonumber \\
&& \quad \times\frac{\Gamma(k+\alpha+\beta+1)}{\Gamma(k+\alpha+1)}
P_{k}^{(\alpha,\beta)}(z).
\label{A.26}
\end{eqnarray}
From it, after exploiting elementary properties of the gamma and the
digamma functions \cite{Davi65,Magn66,Grad94}, one obtains two
related formulas:
\begin{eqnarray}
\frac{\partial P_{n}^{(\alpha,\beta)}(z)}{\partial\beta}
&=& [\psi(2n+\alpha+\beta+1)-\psi(n+\alpha+\beta+1)]
P_{n}^{(\alpha,\beta)}(z)
\nonumber \\
&& +\,\frac{1}{\Gamma(-n-\alpha)\Gamma(n+\alpha+\beta+1)}
\sum_{k=0}^{n-1}
\frac{2k+\alpha+\beta+1}{(n-k)(k+n+\alpha+\beta+1)}
\nonumber \\
&& \quad \times\Gamma(-k-\alpha)\Gamma(k+\alpha+\beta+1)
P_{k}^{(\alpha,\beta)}(z)
\label{A.27}
\end{eqnarray}
and
\begin{eqnarray}
\frac{\partial P_{n}^{(\alpha,\beta)}(z)}{\partial\beta}
&=& [\psi(-2n-\alpha-\beta)-\psi(-n-\alpha-\beta)]
P_{n}^{(\alpha,\beta)}(z)
\nonumber \\
&& +\,(-)^{n}\frac{\Gamma(-n-\alpha-\beta)}{\Gamma(-n-\alpha)}
\sum_{k=0}^{n-1}(-)^{k}
\frac{2k+\alpha+\beta+1}{(n-k)(k+n+\alpha+\beta+1)}
\nonumber \\
&& \quad \times\frac{\Gamma(-k-\alpha)}{\Gamma(-k-\alpha-\beta)}
P_{k}^{(\alpha,\beta)}(z).
\label{A.28}
\end{eqnarray}

The following special cases of Eqs.\ (\ref{A.15}) to (\ref{A.28})
have been exploited by us in sections {\ref{V.1.2} and \ref{V.4.2}
(notice the constraints imposed on $m$ in some of the formulas listed
below):
\begin{eqnarray}
\frac{\partial P_{n}^{(m,\beta)}(z)}{\partial\beta}\bigg|_{\beta=-m}
&=& -\,P_{n}^{(m,-m)}(z)\ln\frac{z+1}{2}
\nonumber \\
&& +\,\frac{1}{2^{n}n!}\left(\frac{z+1}{z-1}\right)^{m}
\frac{\mathrm{d}^{n}}{\mathrm{d}z^{n}}
\left[(z-1)^{n+m}(z+1)^{n-m}\ln\frac{z+1}{2}\right],
\label{A.29}
\end{eqnarray}
\begin{eqnarray}
\frac{\partial P_{n}^{(m,\beta)}(z)}{\partial\beta}\bigg|_{\beta=-m}
&=& -\,\psi(n+1)P_{n}^{(m,-m)}(z)
\nonumber \\
&& +\,\frac{(n+m)!}{n!}\sum_{k=0}^{n}
\frac{(k+n)!\psi(k+n+1)}{k!(k+m)!(n-k)!}
\left(\frac{z-1}{2}\right)^{k},
\label{A.30}
\end{eqnarray}
\begin{eqnarray}
\frac{\partial P_{n}^{(m,\beta)}(z)}{\partial\beta}\bigg|_{\beta=-m}
&=& -\,[\psi(n+1)-\psi(n-m+1)]P_{n}^{(m,-m)}(z)
\nonumber \\
&& -\,(-)^{n+m}\frac{(n-m)!}{n!}\sum_{k=0}^{m-1}
\frac{(k+n)!(m-k-1)!}{k!(n-k)!}\left(\frac{z+1}{2}\right)^{k}
\nonumber \\
&& +\,(-)^{n+m}\frac{(n-m)!}{n!}\left(\frac{z+1}{2}\right)^{m}
\sum_{k=0}^{n-m}(-)^{k}\frac{(k+n+m)!}{k!(k+m)!(n-m-k)!}
\nonumber \\
&& \quad \times[\psi(k+n+m+1)-\psi(k+1)]
\left(\frac{z+1}{2}\right)^{k}
\qquad (0\leqslant m\leqslant n),
\nonumber \\
&&
\label{A.31}
\end{eqnarray}
\begin{eqnarray}
\frac{\partial P_{n}^{(m,\beta)}(z)}
{\partial\beta}\bigg|_{\beta=-m}
&=& [\psi(m-n)-\psi(n+1)]P_{n}^{(m,-m)}(z)
\nonumber \\
&& +\,\frac{1}{n!(m-n-1)!}
\sum_{k=0}^{n}\frac{(k+n)!(m-k-1)!}{k!(n-k)!}
\nonumber \\
&& \quad 
\times[\psi(k+n+1)-\psi(m-k)]\left(\frac{z+1}{2}\right)^{k}
\qquad (m>n),
\label{A.32}
\end{eqnarray}
\begin{eqnarray}
\frac{\partial P_{n}^{(m,\beta)}(z)}
{\partial\beta}\bigg|_{\beta=-m}
&=& \psi(n-m+1)P_{n}^{(m,-m)}(z)
\nonumber \\
&& -\,(n-m)!(n+m)!\left(\frac{z+1}{2}\right)^{n}
\nonumber \\
&& \quad \times\sum_{k=0}^{n-m}
\frac{\psi(n-m-k+1)}{k!(k+m)!(n-k)!(n-m-k)!}
\left(\frac{z-1}{z+1}\right)^{k}
\nonumber \\
&& -\,(n-m)!(n+m)!\left(\frac{z+1}{z-1}\right)^{m}
\left(\frac{z-1}{2}\right)^{n}
\nonumber \\
&& \quad 
\times\sum_{k=1}^{m}(-)^{k}\frac{(k-1)!}{(k+n)!(k+n-m)!(m-k)!}
\left(\frac{z-1}{z+1}\right)^{k}
\qquad (0\leqslant m\leqslant n),
\nonumber \\
&&
\label{A.33}
\end{eqnarray}
\begin{eqnarray}
\frac{\partial P_{n}^{(m,\beta)}(z)}
{\partial\beta}\bigg|_{\beta=-m}
&=& \psi(m-n)P_{n}^{(m,-m)}(z)
\nonumber \\
&& -\,\frac{(n+m)!}{(m-n-1)!}\left(\frac{z+1}{2}\right)^{n}
\nonumber \\
&& \quad \times\sum_{k=0}^{n}(-)^{k}
\frac{(k+m-n-1)!\psi(k+m-n)}{k!(k+m)!(n-k)!}
\left(\frac{z-1}{z+1}\right)^{k}
\qquad (m>n),
\nonumber \\
&&
\label{A.34}
\end{eqnarray}
\begin{eqnarray}
\frac{\partial P_{n}^{(m,\beta)}(z)}
{\partial\beta}\bigg|_{\beta=-m}
&=& [\psi(2n+1)-\psi(n+1)]P_{n}^{(m,-m)}(z)
\nonumber \\
&& +\,(-)^{n}\frac{(n+m)!}{n!}
\sum_{k=0}^{n-1}(-)^{k}\frac{2k+1}{(n-k)(k+n+1)}
\frac{k!}{(k+m)!}P_{k}^{(m,-m)}(z),
\nonumber \\
&&
\label{A.35}
\end{eqnarray}
\begin{eqnarray}
\frac{\partial P_{n}^{(-m,\beta)}(z)}{\partial\beta}\bigg|_{\beta=m}
&=& -\,P_{n}^{(-m,m)}(z)\ln\frac{z+1}{2}
\nonumber \\
&& +\,\frac{1}{2^{n}n!}\left(\frac{z-1}{z+1}\right)^{m}
\frac{\mathrm{d}^{n}}{\mathrm{d}z^{n}}
\left[(z-1)^{n-m}(z+1)^{n+m}\ln\frac{z+1}{2}\right],
\label{A.36}
\end{eqnarray}
\begin{eqnarray}
\frac{\partial P_{n}^{(-m,\beta)}(z)}
{\partial\beta}\bigg|_{\beta=m} &=& -\psi(n+1)P_{n}^{(-m,m)}(z)
\nonumber \\
&& +\,\frac{(n-m)!}{n!}\left(\frac{z-1}{2}\right)^{m}
\sum_{k=0}^{n-m}\frac{(k+n+m)!\psi(k+n+m+1)}{k!(k+m)!(n-m-k)!}
\left(\frac{z-1}{2}\right)^{k}
\nonumber \\
&& (0\leqslant m\leqslant n),
\label{A.37}
\end{eqnarray}
\begin{eqnarray}
\frac{\partial P_{n}^{(-m,\beta)}(z)}{\partial\beta}\bigg|_{\beta=m}
&=& -\,\psi(n+1)P_{n}^{(-m,m)}(z)
\nonumber \\
&& +\,\frac{(-)^{n}}{n!(m-n-1)!}
\sum_{k=0}^{n}(-)^{k}\frac{(k+n)!(m-k-1)!\psi(k+n+1)}{k!(n-k)!}
\left(\frac{z-1}{2}\right)^{k}
\nonumber \\
&& (m>n),
\label{A.38}
\end{eqnarray}
\begin{eqnarray}
\frac{\partial P_{n}^{(-m,\beta)}(z)}{\partial\beta}\bigg|_{\beta=m}
&=& [\psi(n+m+1)-\psi(n+1)]P_{n}^{(-m,m)}(z)
\nonumber \\
&& +\,(-)^{n}\frac{(n+m)!}{n!}\sum_{k=0}^{n}(-)^{k}
\frac{(k+n)!}{k!(k+m)!(n-k)!}
\nonumber \\
&& \quad \times[\psi(k+n+1)-\psi(k+m+1)]
\left(\frac{z+1}{2}\right)^{k},
\label{A.39}
\end{eqnarray}
\begin{eqnarray}
\frac{\partial P_{n}^{(-m,\beta)}(z)}{\partial\beta}\bigg|_{\beta=m}
&=& \psi(n+m+1)P_{n}^{(-m,m)}(z)
\nonumber \\
&& -\,(n-m)!(n+m)!\left(\frac{z-1}{z+1}\right)^{m}
\left(\frac{z+1}{2}\right)^{n}
\nonumber \\
&& \quad \times\sum_{k=0}^{n-m}
\frac{\psi(n-k+1)}{k!(k+m)!(n-k)!(n-m-k)!}
\left(\frac{z-1}{z+1}\right)^{k}
\qquad (0\leqslant m\leqslant n),
\nonumber \\
&&
\label{A.40}
\end{eqnarray}
\begin{eqnarray}
\frac{\partial P_{n}^{(-m,\beta)}(z)}{\partial\beta}\bigg|_{\beta=m}
&=& \psi(n+m+1)P_{n}^{(-m,m)}(z)
\nonumber \\
&& -\,(-)^{n}\frac{(n+m)!}{(m-n-1)!}\left(\frac{z+1}{2}\right)^{n}
\nonumber \\
&& \quad \times\sum_{k=0}^{n}(-)^{k}
\frac{(m-k-1)!\psi(n+m-k+1)}{k!(n-k)!(n+m-k)!}
\left(\frac{z-1}{z+1}\right)^{k}
\qquad (m>n),
\nonumber \\
&&
\label{A.41}
\end{eqnarray}
\begin{eqnarray}
\frac{\partial P_{n}^{(-m,\beta)}(z)}{\partial\beta}\bigg|_{\beta=m}
&=& [\psi(2n+1)-\psi(n+1)]P_{n}^{(-m,m)}(z)
\nonumber \\
&& +\,(-)^{n+m}\frac{(n-m)!}{n!}
\sum_{k=0}^{n-m-1}(-)^{k}\frac{2k+2m+1}{(n-m-k)(k+n+m+1)}
\nonumber \\
&& \quad \times\frac{(k+m)!}{k!}P_{k+m}^{(-m,m)}(z)
\qquad (0\leqslant m\leqslant n),
\label{A.42}
\end{eqnarray}
\begin{eqnarray}
\frac{\partial P_{n}^{(-m,\beta)}(z)}{\partial\beta}\bigg|_{\beta=m}
&=& [\psi(2n+1)-\psi(n+1)]P_{n}^{(-m,m)}(z)
\nonumber \\
&& +\,\frac{1}{n!(m-n-1)!}
\sum_{k=0}^{n-1}\frac{2k+1}{(n-k)(k+n+1)}
\nonumber \\
&& \quad \times k!(m-k-1)!P_{k}^{(-m,m)}(z)
\qquad (m>n),
\label{A.43}
\end{eqnarray}
\begin{eqnarray}
\frac{\partial P_{n-m}^{(m,\beta)}(z)}{\partial\beta}\bigg|_{\beta=m}
&=& -\,P_{n-m}^{(m,m)}(z)\ln\frac{z+1}{2}
\nonumber \\
&& +\,\frac{1}{2^{n-m}(n-m)!}(z^{2}-1)^{-m}
\frac{\mathrm{d}^{n-m}}{\mathrm{d}z^{n-m}}
\left[(z^{2}-1)^{n}\ln\frac{z+1}{2}\right]
\nonumber \\
&& (0\leqslant m\leqslant n),
\label{A.44}
\end{eqnarray}
\begin{eqnarray}
\frac{\partial P_{n-m}^{(m,\beta)}(z)}
{\partial\beta}\bigg|_{\beta=m}
&=& -\,\psi(n+m+1)P_{n-m}^{(m,m)}(z)
\nonumber \\
&& +\,\frac{n!}{(n+m)!}
\sum_{k=0}^{n-m}\frac{(k+n+m)!\psi(k+n+m+1)}{k!(k+m)!(n-m-k)!}
\left(\frac{z-1}{2}\right)^{k}
\nonumber \\
&& (0\leqslant m\leqslant n),
\label{A.45}
\end{eqnarray}
\begin{eqnarray}
\frac{\partial P_{n-m}^{(m,\beta)}(z)}
{\partial\beta}\bigg|_{\beta=m}
&=& -\,[\psi(n+m+1)-\psi(n+1)]P_{n-m}^{(m,m)}(z)
\nonumber \\
&& +\,(-)^{n+m}\frac{n!}{(n+m)!}
\sum_{k=0}^{n-m}(-)^{k}\frac{(k+n+m)!}{k!(k+m)!(n-m-k)!}
\nonumber \\
&& \quad \times[\psi(k+n+m+1)-\psi(k+m+1)]
\left(\frac{z+1}{2}\right)^{k}
\qquad (0\leqslant m\leqslant n),
\nonumber \\
&&
\label{A.46}
\end{eqnarray}
\begin{eqnarray}
\frac{\partial P_{n-m}^{(m,\beta)}(z)}
{\partial\beta}\bigg|_{\beta=m}
&=& \psi(n+1)P_{n-m}^{(m,m)}(z)
\nonumber \\
&& -\,(n!)^{2}\left(\frac{z+1}{2}\right)^{n-m}
\sum_{k=0}^{n-m}\frac{\psi(n-k+1)}{k!(k+m)!(n-k)!(n-m-k)!}
\left(\frac{z-1}{z+1}\right)^{k}
\nonumber \\
&& (0\leqslant m\leqslant n),
\label{A.47}
\end{eqnarray}
\begin{eqnarray}
\frac{\partial P_{n+m}^{(-m,\beta)}(z)}
{\partial\beta}\bigg|_{\beta=-m}
&=& -\,P_{n+m}^{(-m,-m)}(z)\ln\frac{z+1}{2}
\nonumber \\
&& +\,\frac{1}{2^{n+m}(n+m)!}(z^{2}-1)^{m}
\frac{\mathrm{d}^{n+m}}{\mathrm{d}z^{n+m}}
\left[(z^{2}-1)^{n}\ln\frac{z+1}{2}\right],
\label{A.48}
\end{eqnarray}
\begin{eqnarray}
\frac{\partial P_{n+m}^{(-m,\beta)}(z)}
{\partial\beta}\bigg|_{\beta=-m}
&=& -\,\psi(n-m+1)P_{n+m}^{(-m,-m)}(z)
\nonumber \\
&& +\,\frac{n!}{(n-m)!}\left(\frac{z-1}{2}\right)^{m}\sum_{k=0}^{n}
\frac{(k+n)!\psi(k+n+1)}{k!(k+m)!(n-k)!}
\left(\frac{z-1}{2}\right)^{k}
\nonumber \\
&& (0\leqslant m\leqslant n),
\label{A.49}
\end{eqnarray}
\begin{eqnarray}
\frac{\partial P_{n+m}^{(-m,\beta)}(z)}
{\partial\beta}\bigg|_{\beta=-m}
&=& [\psi(n+1)-\psi(n-m+1)]P_{n+m}^{(-m,-m)}(z)
\nonumber \\
&& -\,(-)^{n}\frac{n!}{(n-m)!}
\sum_{k=0}^{m-1}\frac{(k+n-m)!(m-k-1)!}{k!(n+m-k)!}
\left(\frac{z+1}{2}\right)^{k}
\nonumber \\
&& +\,(-)^{n}\frac{n!}{(n-m)!}\left(\frac{z+1}{2}\right)^{m}
\sum_{k=0}^{n}(-)^{k}\frac{(k+n)!}{k!(k+m)!(n-k)!}
\nonumber \\
&& \quad \times[\psi(k+n+1)-\psi(k+1)]
\left(\frac{z+1}{2}\right)^{k}
\qquad (0\leqslant m\leqslant n),
\nonumber \\
&&
\label{A.50}
\end{eqnarray}
\begin{eqnarray}
\frac{\partial P_{n+m}^{(-m,\beta)}(z)}
{\partial\beta}\bigg|_{\beta=-m}
&=& \psi(n+1)P_{n+m}^{(-m,-m)}(z)
\nonumber \\
&& -\,(n!)^{2}\left(\frac{z-1}{2}\right)^{m}
\left(\frac{z+1}{2}\right)^{n}
\nonumber \\
&& \quad \times\sum_{k=0}^{n-m}
\frac{\psi(n-m-k+1)}{k!(k+m)!(n-k)!(n-m-k)!}
\left(\frac{z-1}{z+1}\right)^{k}
\nonumber \\
&& -\,(n!)^{2}\left(\frac{z-1}{2}\right)^{n}
\left(\frac{z+1}{2}\right)^{m}
\nonumber \\
&& \quad \times\sum_{k=1}^{m}(-)^{k}
\frac{(k-1)!}{(k+n)!(k+n-m)!(m-k)!}\left(\frac{z-1}{z+1}\right)^{k}
\qquad (0\leqslant m\leqslant n),
\nonumber \\
&&
\label{A.51}
\end{eqnarray}
\begin{eqnarray}
\frac{\partial P_{m-n-1}^{(-m,\beta)}(z)}
{\partial\beta}\bigg|_{\beta=-m}
&=& -\,P_{m-n-1}^{(-m,-m)}(z)\ln\frac{z+1}{2}
\nonumber \\
&& +\,\frac{1}{2^{m-n-1}(m-n-1)!}(z^{2}-1)^{m}
\frac{\mathrm{d}^{m-n-1}}{\mathrm{d}z^{m-n-1}}
\left[(z^{2}-1)^{-n-1}\ln\frac{z+1}{2}\right]
\nonumber \\
&& (m>n),
\label{A.52}
\end{eqnarray}
\begin{eqnarray}
\frac{\partial P_{m-n-1}^{(-m,\beta)}(z)}
{\partial\beta}\bigg|_{\beta=-m}
&=& -\,\psi(n+m+1)P_{m-n-1}^{(-m,-m)}(z)
\nonumber \\
&& -\,(-)^{n+m}\frac{(n+m)!}{n!}
\sum_{k=0}^{m-n-1}\frac{(m-k-1)!\psi(n+m-k+1)}{k!(n+m-k)!(m-n-k-1)!}
\left(\frac{z-1}{2}\right)^{k}
\nonumber \\
&& (m>n),
\label{A.53}
\end{eqnarray}
\begin{eqnarray}
\frac{\partial P_{m-n-1}^{(-m,\beta)}(z)}
{\partial\beta}\bigg|_{\beta=-m}
&=& -\,[\psi(n+m+1)-\psi(n+1)]P_{m-n-1}^{(-m,-m)}(z)
\nonumber \\
&& +\,\frac{(n+m)!}{n!}
\sum_{k=0}^{m-n-1}(-)^{k}\frac{(m-k-1)!}
{k!(n+m-k)!(m-n-k-1)!}
\nonumber \\
&& \quad \times[\psi(n+m-k+1)-\psi(m-k)]
\left(\frac{z+1}{2}\right)^{k}
\qquad (m>n),
\nonumber \\
&&
\label{A.54}
\end{eqnarray}
\begin{eqnarray}
\frac{\partial P_{m-n-1}^{(-m,\beta)}(z)}
{\partial\beta}\bigg|_{\beta=-m}
&=& \psi(n+1)P_{m-n-1}^{(-m,-m)}(z)
\nonumber \\
&& +\,\frac{(-)^{n+m}}{(n!)^{2}}\left(\frac{z+1}{2}\right)^{m-n-1}
\nonumber \\
&& \quad \times\sum_{k=0}^{m-n-1}
\frac{(k+n)!(m-k-1)!\psi(k+n+1)}{k!(m-n-k-1)!}
\left(\frac{z-1}{z+1}\right)^{k}
\qquad (m>n),
\nonumber \\
&&
\label{A.55}
\end{eqnarray}
\begin{eqnarray}
\frac{\partial P_{m-n-1}^{(-m,\beta)}(z)}
{\partial\beta}\bigg|_{\beta=-m}
&=& -\,[\psi(n+m+1)-\psi(2n+2)]P_{m-n-1}^{(-m,-m)}(z)
\nonumber \\
&& -\,(-)^{n+m}\frac{(n+m)!}{n!}
\sum_{k=0}^{m-n-2}(-)^{k}\frac{2m-2k-1}{(n+m-k)(m-n-k-1)}
\nonumber \\
&& \quad \times\frac{(m-k-1)!}{(2m-k-1)!}P_{k}^{(-m,-m)}(z)
\qquad (m>n).
\label{A.56}
\end{eqnarray}
%
%

%

\begin{thebibliography}{99}
\bibitem{Szmy06a}
   R.\ Szmytkowski,
   \textit{On the derivative of the Legendre function of the first 
   kind with respect to its degree},
   J.\ Phys.\ A 39 (2006) 15147;
   corrigendum: 40 (2007) 7819.
\bibitem{Joll19}
   A.\ E.\ Jolliffe,
   \textit{A form for $\frac{\mathrm{d}}{\mathrm{d}n}P_{n}(\mu)$, 
   where $P_{n}(\mu)$ is the Legendre polynomial of degree $n$},
   Mess.\ Math.\ 49 (1919) 125.
\bibitem{Brom13}
   T.\ J.\ I'A Bromwich,
   \textit{Certain potential functions and a new solution of 
   Laplace's equation},
   Proc.\ Lond.\ Math.\ Soc.\ 12 (1913) 100.
\bibitem{Sche41}
   S.\ A.\ Schelkunoff,
   \textit{Theory of antennas of arbitrary size and shape},
   Proc.\ IRE 29 (1941) 493;
   corrigendum: 31 (1943) 38;
   reprint: Proc.\ IEEE 72 (1984) 1165.
   Notice that Schelkunoff used a definition of the digamma function
   different from that adopted in \cite{Szmy06a,Szmy07} and in the 
   present paper.
\bibitem{Davi65}
   P.\ J.\ Davis,
   \textit{Gamma function and related functions\/},
   in \textit{Handbook of Mathematical Functions},
   edited by M.\ Abramowitz and I.\ A.\ Stegun,
   Dover, New York, p.\ 253.
\bibitem{Magn66}
   W.\ Magnus, F.\ Oberhettinger, and R.\ P.\ Soni,
   \textit{Formulas and Theorems for the Special Functions of
   Mathematical Physics\/} 3rd edn,
   Springer, Berlin, 1966.
\bibitem{Grad94}
   I.\ S.\ Gradshteyn and I.\ M.\ Ryzhik,
   \emph{Table of Integrals, Series, and Products\/} 5th edn,
   Academic, San Diego, 1994.
   In section 8.76, a few formulas for the derivative of the
   associated Legendre function of the first kind with respect to its
   degree (and \emph{not\/} order, as incorrectly the title of the
   section says) are listed.
\bibitem{Szmy07}
   R.\ Szmytkowski,
   \textit{Addendum to `On the derivative of the Legendre function 
   of the first kind with respect to its degree'}, 
   J.\ Phys.\ A 40 (2007) 14887.  We note parenthetically
   that the simplest way of proving Eq.\ (5) in that paper is to
   differentiate both sides of the identity
   $\mathrm{d}^{k}z^{\alpha}/\mathrm{d}z^{k}
   =[\Gamma(\alpha+1)/\Gamma(\alpha-k+1)]z^{\alpha-k}$ with respect
   to $\alpha$.
\bibitem{Cars13}
   H.\ S.\ Carslaw,
   \textit{Integral equations and the determination of Green's 
   functions in the theory of potential},
   Proc.\ Edinburgh Math.\ Soc.\ 31 (1913) 71.
\bibitem{Cars14a}
   H.\ S.\ Carslaw,
   \textit{The scattering of sound waves by a cone},
   Math.\ Ann.\ 75 (1914) 133;
   corrigendum: 75 (1914) 592.
\bibitem{Cars14b}
   H.\ S.\ Carslaw,
   \textit{The Green's function for the equation
   $\boldsymbol{\nabla}^{2}u+k^{2}u=0$},
   Proc.\ Lond.\ Math.\ Soc.\ 13 (1914) 236.
\bibitem{Macd15}
   H.\ M.\ Macdonald,
   \textit{A class of diffraction problems},
   Proc.\ Lond.\ Math.\ Soc.\ 14 (1915) 410.
\bibitem{Cars21}
   H.\ S.\ Carslaw,
   \textit{Introduction to the Mathematical Theory of the Conduction 
   of Heat in Solids},
   Macmillan, London, 1921, pp.\ 145--147.
\bibitem{Cars47}
   H.\ S.\ Carslaw and J.\ C.\ Jaeger,
   \textit{Conduction of Heat in Solids},
   Clarendon, Oxford, 1947, pp.\ 214 and 318.
\bibitem{Smyt50}
   W.\ R.\ Smythe,
   \textit{Static and Dynamic Electricity\/} 2nd edn,
   McGraw-Hill, New York, 1950, pp.\ 156--157. \\
   W.\ R.\ Smythe,
   \textit{Static and Dynamic Electricity\/} 3rd edn,
   McGraw-Hill, New York, 1968, pp.\ 166--167.
\bibitem{Fels55}
   L.\ B.\ Felsen,
   \textit{Backscattering from wide-angle and narrow-angle cones},
   J.\ Appl.\ Phys.\ 26 (1955) 138.
\bibitem{Bail56}
   L.\ L.\ Bailin and S.\ Silver,
   \textit{Exterior electromagnetic boundary value problems for 
   spheres and cones},
   IRE Trans.\ Antennas Propag.\ 4 (1956) 5;
   corrigendum: 5 (1957) 313.
\bibitem{Fels57}
   L.\ B.\ Felsen,
   \textit{Plane-wave scattering by small-angle cones},
   IRE Trans.\ Antennas Propag.\ 5 (1957) 121.
\bibitem{Fels59}
   L.\ B.\ Felsen,
   \textit{Radiation from ring sources in the presence of a 
   semi-infinite cone},
   IRE Trans.\ Antennas Propag.\ 7 (1959) 168;
   corrigendum: 7 (1959) 251.
\bibitem{Jone64}
   D.\ S.\ Jones,
   \textit{The Theory of Electromagnetism},
   Pergamon, Oxford, 1964, p.\ 614.
\bibitem{Bowm69}
   J.\ J.\ Bowman,
   \textit{The cone\/},
   in \textit{Electromagnetic and Acoustic Scattering by Simple
   Shapes},
   edited by J.\ J.\ Bowman, T.\ B.\ A.\ Senior and 
   P.\ L.\ E.\ Uslenghi,
   North-Holland, Amsterdam, 1969, p.\ 637.
\bibitem{Fels73}
   L.\ B.\ Felsen and N.\ Marcuvitz,
   \textit{Radiation and Scattering of Waves},
   Prentice-Hall, Englewood Cliffs, NJ, 1973 
   (reprint: IEEE Press, Piscataway, NJ, 1994),
   pp.\ 320, 321, 703 and 734.
\bibitem{Ariy85}
   J.\ C.\ Ariyasu and D.\ L.\ Mills,
   \textit{Inelastic electron scattering by long-wavelength, acoustic
   phonons; image potential modulation as a mechanism},
   Surf.\ Sci.\ 155 (1985) 607, appendix~B.
\bibitem{Baue87}
   H.\ F.\ Bauer,
   \textit{Mass transport in a three-dimensional diffusor or
   confusor},
   W{\"a}rme-Stoff{\"u}bertrag.\ 21 (1987) 51.
\bibitem{Jone86}
   D.\ S.\ Jones,
   \textit{Acoustic and Electromagnetic Waves},
   Clarendon, Oxford, 1986, p.\ 591.
\bibitem{Baue92}
   H.\ F.\ Bauer,
   \textit{Response of axially excited spherical and conical liquid 
   systems with anchored edges},
   Forsch.\ Ing.-Wes.\ 58 (4) (1992) 96.
\bibitem{Broa99}
   E.\ G.\ Broadbent and D.\ W.\ Moore,
   \textit{The inclination of a hollow vortex with an inclined plane 
   and the acoustic radiation produced},
   Proc.\ R.\ Soc.\ Lond.\ A 455 (1999) 1979.
\bibitem{VanB07}
   J.\ Van Bladel,
   \emph{Electromagnetic Fields\/} 2nd edn,
   IEEE Press, Piscataway, 2007, section 16.7.1.
\bibitem{Todh75}
   I.\ Todhunter,
   \textit{An Elementary Treatise on Laplace's Functions, Lam{\'e}'s
   Functions and Bessel's Functions},
   Macmillan, London, 1875.
\bibitem{Ferr77}
   N.\ M.\ Ferrers,
   \textit{An Elementary Treatise on Spherical Harmonics},
   Macmillan, London, 1877.
\bibitem{Neum78}
   F.\ Neumann,
   \textit{Beitr{\"a}ge zur Theorie der Kugelfunctionen},
   Teubner, Leipzig, 1878.
\bibitem{Hein78}
   E.\ Heine,
   \textit{Handbuch der Kugelfunctionen\/} vol.~1, 2nd edn,
   Reimer, Berlin, 1878.
\bibitem{Hein81}
   E.\ Heine,
   \textit{Handbuch der Kugelfunctionen\/} vol.~2, 2nd edn,
   Reimer, Berlin, 1881.
\bibitem{Olbr87}
   R.\ Olbricht,
   \textit{Studien {\"u}ber die Kugel- und Cylinderfunctionen},
   Nova Acta Leop.\ Carol.\ Akad.\ 52 (1887) 1.
\bibitem{Byer93}
   W.\ E.\ Byerly,
   \textit{An Elementary Treatise on Fourier's Series and Spherical,
   Cylindrical, and Ellipsoidal Harmonics, with Applications to
   Problems in Mathematical Physics},
   Ginn, Boston, 1893 
   (reprint: Dover, Mineola, NY, 2003).
\bibitem{Hobs96}
   E.\ W.\ Hobson,
   \textit{On a type of spherical harmonics of unrestricted degree, 
   order, and argument},
   Phil.\ Trans.\ R.\ Soc.\ Lond.\ A 187 (1896) 443.
\bibitem{Wang04}
   A.\ Wangerin,
   \textit{Theorie der Kugelfunktionen und der verwandten Funktionen,
   insbesondere der Lam{\'e}'schen und Bessel'schen\/},
   in \textit{Encyklop{\"a}die der mathematischen Wissenschaften\/}
   vol.~2.1, Teubner, Leipzig, 1904, p.\ 695.
\bibitem{Barn07}
   E.\ W.\ Barnes,
   \textit{On generalized Legendre functions},
   Quart.\ J.\ Pure Appl.\ Math.\ 39 (1907) 97.
   The associated Legendre function of the second kind defined in
   that work differs from the counterpart function of
   Hobson \cite{Hobs96,Hobs31} used in the present paper. The
   relationship between the two functions is:
   $[Q_{\nu}^{\mu}(z)]_{\mathrm{Barnes}}
   =\{\mathrm{e}^{-\mathrm{i}\pi\mu}\sin[\pi(\nu+\mu)]/\sin(\pi\nu)\}
   Q_{\nu}^{\mu}(z)$.
\bibitem{Wang21}
   A.\ Wangerin,
   \textit{Theorie des Potentials und der Kugelfunktionen\/} vol.~2,
   de Gruyter, Berlin, 1921.
\bibitem{Hobs31}
   E.\ W.\ Hobson,
   \textit{The Theory of Spherical and Ellipsoidal Harmonics},
   Cambridge University Press, Cambridge, 1931
   (reprint: Chelsea, New York, 1955).
\bibitem{Pras30}
   G.\ Prasad,
   \textit{A Treatise on Spherical Harmonics and the Functions of
   Bessel and Lam{\'e}, Part I (Elementary)\/},
   Benares Mathematical Society, Benares, 1930.
\bibitem{Pras32}
   G.\ Prasad,
   \textit{A Treatise on Spherical Harmonics and the Functions of
   Bessel and Lam{\'e}, Part II (Advanced)\/},
   Benares Mathematical Society, Benares, 1932.
\bibitem{Snow52}
   Ch.\ Snow,
   \textit{Hypergeometric and Legendre Functions with Applications to
   Integral Equations of Potential Theory\/} 2nd edn,
   National Bureau of Standards, Washington, DC, 1952.
   The associated Legendre functions defined in that book differ 
   from the counterpart functions of Hobson \cite{Hobs96,Hobs31} 
   used in the present paper. The relationships between the two sets 
   of functions are: $[P_{\nu}^{\mu}(z)]_{\mathrm{Snow}}
   =[\Gamma(\nu+\mu+1)/\Gamma(\nu-\mu+1)]P_{\nu}^{-\mu}(z)$ and
   $[Q_{\nu}^{\mu}(z)]_{\mathrm{Snow}}
   =\mathrm{e}^{-\mathrm{i}\pi\mu}\cos(\pi\mu)Q_{\nu}^{\mu}(z)$.
\bibitem{Lens54}
   J.\ Lense,
   \textit{Kugelfunktionen\/} 2nd edn,
   Geest \& Portig, Leipzig, 1954.
\bibitem{Robi57}
   L.\ Robin,
   \textit{Fonctions Sph{\'e}riques de Legendre et Fonctions
   Sph{\'e}ro{\"{\i}}dales\/} vol.~1,
   Gauthier-Villars, Paris, 1957.
\bibitem{Robi58}
   L.\ Robin,
   \textit{Fonctions Sph{\'e}riques de Legendre et Fonctions
   Sph{\'e}ro{\"{\i}}dales\/} vol.~2,
   Gauthier-Villars, Paris, 1958.
\bibitem{Robi59}
   L.\ Robin,
   \textit{Fonctions Sph{\'e}riques de Legendre et Fonctions
   Sph{\'e}ro{\"{\i}}dales\/} vol.~3,
   Gauthier-Villars, Paris, 1959.
\bibitem{Krat60}
   A.\ Kratzer and W.\ Franz,
   \textit{Transzendente Funktionen},
   Akademische Verlagsgesellschaft, Leipzig, 1960, chapter 5.
\bibitem{MacR67}
   T.\ M.\ MacRobert,
   \textit{Spherical Harmonics\/} 3rd edn,
   Pergamon, Oxford, 1967.
\bibitem{Wang89}
   Z.\ X.\ Wang and D.\ R.\ Guo,
   \textit{Special Functions},
   World Scientific, Singapore, 1989, chapter 5.
\bibitem{Temm96}
   N.\ M.\ Temme,
   \textit{Special Functions.\ An Introduction to the Classical
   Functions of Mathematical Physics},
   Wiley, New York, 1996, chapter 8.
\bibitem{Magn48}
   W.\ Magnus and F.\ Oberhettinger,
   \textit{Formeln und S{\"a}tze f{\"u}r die speziellen Funktionen 
   der mathematischen Physik\/} 2nd edn,
   Springer, Berlin, 1948.
\bibitem{Erde53}
   A.\ Erd{\'e}lyi (ed.),
   \textit{Higher Transcendental Functions\/} vol.~1,
   McGraw-Hill, New York, 1953, chapter III.
\bibitem{Jahn60}
   E.\ Jahnke, F.\ Emde, and F.\ L{\"o}sch,
   \textit{Tafeln h{\"o}herer Funktionen\/} 6th edn,
   Teubner, Stuttgart, 1960.
\bibitem{Steg65}
   I.\ A.\ Stegun,
   \textit{Legendre functions\/},
   in \textit{Handbook of Mathematical Functions},
   edited by M.\ Abramowitz and I.\ A.\ Stegun,
   Dover, New York, 1965, p.\ 331.
\bibitem{Prud83}
   A.\ P.\ Prudnikov, Yu.\ A.\ Brychkov, and O.\ I.\ Marichev,
   \textit{Integrals and Series. Special Functions},
   Nauka, Moscow, 1983 (in Russian).
\bibitem{Prud03}
   A.\ P.\ Prudnikov, Yu.\ A.\ Brychkov, and O.\ I.\ Marichev,
   \textit{Integrals and Series. Special Functions. Supplementary
   Chapters\/} 2nd edn,
   Fizmatlit, Moscow, 2003 (in Russian).
\bibitem{Bryc08}
   Yu.\ A.\ Brychkov,
   \textit{Handbook of Special Functions. Derivatives, Integrals,
   Series and Other Formulas},
   Chapman \& Hall/CRC, Boca Raton, FL, 2008.
   The section on derivatives of the associated Legendre function of
   the first kind with respect to its parameters was absent in an
   earlier Russian edition of this book (Fizmatlit, Moscow, 2006).
\bibitem{Szmy09a}
   R.\ Szmytkowski,
   \textit{On the derivative of the associated Legendre function of 
   the first kind of integer degree with respect to its order (with
   applications to the construction of the associated Legendre
   function of the second kind of integer degree and order)},
   J.\ Math.\ Chem.\ 46 (2009) 231.
\bibitem{Macd00}
   H.\ M.\ Macdonald,
   \textit{Demonstration of Green's formula for electric density near
   the vertex of a right cone},
   Trans.\ Camb.\ Phil.\ Soc.\ 18 (1900) 292.
\bibitem{Macd02}
   H.\ M.\ Macdonald,
   \textit{Electric Waves},
   Cambridge University Press, Cambridge, 1902, pp.\ 89, 92 and
   95--97.
\bibitem{Lebe39}
   N.\ N.\ Lebedev and M.\ I.\ Kontorovich,
   \textit{On application of inversion formulas to solving some 
   problems of electrodynamics},
   Zh.\ Eksp.\ Teor.\ Fiz.\ 9 (1939) 729 (in Russian).
\bibitem{Grin48}
   G.\ A.\ Grinberg,
   \textit{Selected Problems of Mathematical Theory of Electric and
   Magnetic Phenomena},
   Academy of Sciences of the USSR, Moscow, 1948 (in Russian),
   pp.\ 194--6, 203--4 and 207.
\bibitem{Jean60}
   J.\ Jeans,
   \textit{The Mathematical Theory of Electricity and Magnetism\/}
   5th edn,
   Cambridge University Press, Cambridge, 1960, p.\ 293.
\bibitem{Sakm70}
   I.\ A.\ Sakmar,
   \textit{Direct-channel fermion Regge poles},
   Nuovo Cimento A 70 (1970) 137.
\bibitem{Nort71}
   F.\ H.\ Northover,
   \textit{Applied Diffraction Theory},
   American Elsevier, New York, 1971,
   chapters 8 and 9, and appendices thereto.
\bibitem{vdPo87}
   B.\ van der Pol and H.\ Bremmer,
   \textit{Operational Calculus Based on the Two-sided Laplace
   Integral\/} 3rd edn,
   Chelsea, New York, 1987, p.\ 242.
\bibitem{Soki97}
   A.\ V.\ Sokirko and K.\ B.\ Oldham,
   \textit{The voltammetric response of a conical electrode},
   J.\ Electroanal.\ Chem.\ 430 (1997) 15.
\bibitem{From07}
   P.\ O.\ Fr{\"o}man and S.\ Yngve,
   \textit{A physically important class of integrals expressed as a
   parameter derivative},
   Ann.\ Phys.\ (N.Y.) 322 (2007) 2145.
\bibitem{Coff08}
   M.\ W.\ Coffey,
   \textit{On harmonic binomial series},
   preprint: arXiv:0812.1766v1.
\bibitem{Robi56}
   L.\ Robin,
   \textit{Deriv{\'e}e de la fonction associ{\'e}e de Legendre, de
   premi{\`e}re esp{\`e}ce, par rapport {\`a} son degr{\'e}},
   Compt.\ Rend.\ Acad.\ Sci.\ Paris\ 242 (1956) 57.
\bibitem{Host64}
   L.\ Hostler,
   \textit{Nonrelativistic Coulomb Green's function in momentum
   space},
   J.\ Math.\ Phys.\ 5 (1964) 1235.
\bibitem{Tsu61}
   R.\ Tsu,
   \textit{The evaluation of incomplete normalization integrals and
   derivatives with respect to the order of associated Legendre
   polynomials},
   J.\ Math.\ and Phys.\ 40 (1961) 232.
   The definition of $P_{n}^{m}(x)$ adopted in that paper differs
   from that of Hobson \cite{Hobs31} by the factor $(-)^{m}$. 
   Moreover, what the author called an \emph{order\/} of the Legendre
   function, in the mathematical literature is most commonly named a
   \emph{degree\/} of the latter. Finally, Eq.\ (40) in that paper
   has been misprinted; the innermost differentiation should be with
   respect to $\theta$, \emph{not\/} $\nu$.
\bibitem{Carl87}
   B.\ C.\ Carlson,
   \textit{Dirichlet averages of $x^{t}\log x$},
   SIAM J.\ Math.\ Anal.\ 18 (1987) 550.
\bibitem{Sche77}
   L.\ Schendel,
   \textit{Zusatz zu der Abhandlung {\"u}ber Kugelfunctionen S.\ 86 
   des 80.\ Bandes},
   J.\ Reine Angew.\ Math.\ (Borchardt J.) 82 (1877) 158.
\bibitem{Szmy06b}
   R.\ Szmytkowski,
   \textit{Closed form of the generalized Green's function for the 
   Helmholtz operator on the two-dimensional unit sphere},
   J.\ Math.\ Phys.\ 47 (2006) 063506.
\bibitem{Szeg39}
   G.\ Szeg{\"o},
   \textit{Orthogonal Polynomials},
   American Mathematical Society, New York, 1939, chapter 4.
\bibitem{Froh94}
   J.\ Fr{\"o}hlich,
   \textit{Parameter derivatives of the Jacobi polynomials and the
   Gaussian hypergeometric function},
   Integral Transforms Spec.\ Funct.\ 2 (1994) 253.
\bibitem{Szmy09b}
   R.\ Szmytkowski,
   \textit{A note on parameter derivatives of classical orthogonal
   polynomials},
   preprint: arXiv:0901.2639v1.
\end{thebibliography}
\end{document}